\documentclass[reqno,11pt]{amsart}
\usepackage{amscd,amssymb,verbatim}
\numberwithin{equation}{section}
\theoremstyle{plain}

\newcommand\alp{\alpha}         
\newcommand\bet{\beta}
\newcommand\gam{\gamma}         \newcommand\Gam{\Gamma}
\newcommand\del{\delta}         \newcommand\Del{\Delta}
\newcommand\eps{\varepsilon}
\newcommand\zet{\zeta}

\newcommand\lam{\lambda}                \newcommand\Lam{\Lambda}
\newcommand\sig{\sigma}         \newcommand\Sig{\Sigma}

\newcommand\ome{\omega}         \newcommand\Ome{\Omega}

\newcommand\calE{{\mathcal{E}}}
\newcommand\calF{{\mathcal{F}}}

\newcommand\calH{{\mathcal{H}}}

\newcommand\calL{{\mathcal{L}}}
\newcommand\calM{{\mathcal{M}}}
\newcommand\calN{{\mathcal{N}}}

\newcommand\calP{{\mathcal{P}}}
\newcommand\calQ{{\mathcal{Q}}}
\newcommand\calR{{\mathcal{R}}}
\newcommand\calS{{\mathcal{S}}}
\newcommand\calT{{\mathcal{T}}}
\newcommand\calU{{\mathcal{U}}}

\newcommand\calW{{\mathcal{W}}}

            \newcommand\bfC{{\mathbf C}}

\newcommand\RR{\mathbb{R}}

\newcommand\CC{\mathbb{C}}

\newcommand\nek{,\ldots,}
\newcommand\sdp{\times \hskip -0.3em {\raise 0.3ex
\hbox{$\scriptscriptstyle |$}}} 

\newcommand\Cone{\operatorname{Cone}}

\newcommand\Dom{\operatorname{Dom}}

\newcommand\IM{\operatorname{Im}}

\newcommand\ind{\operatorname{ind}}

\newcommand\Ker{\operatorname{Ker}}

\newcommand\rk{\operatorname{rk}}

\newcommand\supp{\operatorname{supp}}

\newcommand\oJ{{\overline{J}}}

\newcommand\oM{{\overline{M}}}

\newcommand\os{{\bar{s}}}

\newcommand\uc{{\underline{c}}}

\newcommand\tilF{{\widetilde{F}}}

\newcommand\tilh{{\widetilde{h}}}

\newcommand\tilm{{\widetilde{m}}}
\newcommand\tilM{{\widetilde{M}}}

\newcommand\tilN{{\widetilde{N}}}

\newcommand\tilU{{\widetilde{U}}}

\newcommand\tilV{{\widetilde{V}}}

\newcommand\tilW{{\widetilde{W}}}

\newcommand\tilome{{\widetilde{\ome}}}

\renewcommand\tilome{{\widetilde{\ome}}}
\newcommand\tilOme{{\widetilde{\Ome}}}

\renewcommand{\>}{\rangle}
\newcommand{\<}{\langle}

\theoremstyle{plain}
\newtheorem{Thm}[subsection]{Theorem}
\newtheorem{Cor}[subsection]{Corollary}
\newtheorem{Lem}[subsection]{Lemma}
\newtheorem{Prop}[subsection]{Proposition}
\newtheorem{Conjec}[subsection]{Conjecture}

\newtheorem{Def}[subsection]{Definition}

\theoremstyle{remark}

\newtheorem{Rem}[subsection]{Remark}

\def\TeXref#1{%
        \leavevmode\vadjust{\setbox0=\hbox{{\tt
                \  {\tiny \textrm #1}}}%
        \theight=\ht0
        \advance\theight by \lineskip
        \kern -\theight \vbox to
        \theight{\rightline{\rlap{\box0}}%
        \vss}%
        }}%
\newif\ifShowLabels
\ShowLabelstrue
\newdimen\theight
\def\TeXrefEq#1{%
        \leavevmode\vadjust{\setbox0=\hbox{{\tt
                \  {\tiny \textrm #1}}}%
        \theight=\ht1
        \advance\theight by \lineskip
        \kern -\theight \vbox to
        \theight{\rightline{\rlap{\box0}}%
        \vss}%
        }}%

\newcommand{\refs}[1]{Section ~\ref{S:#1}}
\newcommand{\refss}[1]{Subsection ~\ref{SS:#1}}

\newcommand{\reft}[1]{Theorem ~\ref{T:#1}}
\newcommand{\refl}[1]{Lemma ~\ref{L:#1}}
\newcommand{\refp}[1]{Proposition ~\ref{P:#1}}
\newcommand{\refc}[1]{Corollary ~\ref{C:#1}}

\newcommand{\refd}[1]{Definition ~\ref{D:#1}}
\newcommand{\refr}[1]{Remark ~\ref{R:#1}}
\newcommand{\refe}[1]{\eqref{E:#1}}

\newenvironment{thm}[1]%
        { \begin{Thm} \label{T:#1}  \ifShowLabels \TeXref{T:#1} \fi }%
        { \end{Thm} }

\renewcommand{\th}[1]{\begin{thm}{#1}  }
\renewcommand{\eth}{\end{thm} }

\newenvironment{lemma}[1]%
        { \begin{Lem} \label{L:#1}  \ifShowLabels \TeXref{L:#1} \fi }%
        { \end{Lem} }

\newcommand{\lem}[1]{\begin{lemma}{#1} }
\newcommand{\elem}{\end{lemma}}

\newenvironment{propos}[1]%
        { \begin{Prop} \label{P:#1}  \ifShowLabels \TeXref{P:#1} \fi }%
        { \end{Prop} }

\newcommand{\prop}[1]{\begin{propos}{#1} }
\newcommand{\eprop}{\end{propos}}

\newenvironment{corol}[1]%
        { \begin{Cor} \label{C:#1}  \ifShowLabels \TeXref{C:#1} \fi }%
        { \end{Cor} }
\newcommand{\cor}[1]{\begin{corol}{#1}  }
\newcommand{\ecor}{\end{corol}}

\newenvironment{conjec}[1]%
        { \begin{Conjec} \label{Conj:#1}  \ifShowLabels \TeXref{C:#1} \fi }%
        { \end{Conjec} }
\newcommand{\conj}[1]{\begin{conjec}{#1}  }
\newcommand{\econj}{\end{conjec}}

\newenvironment{defeni}[1]%
        { \begin{Def} \label{D:#1}  \ifShowLabels \TeXref{D:#1} \fi }%
        { \end{Def} }
\newcommand{\defe}[1]{\begin{defeni}{#1}  }
\newcommand{\edefe}{\end{defeni}}

\newenvironment{remark}[1]%
        { \begin{Rem} \label{R:#1}  \ifShowLabels \TeXref{R:#1} \fi }%
        { \end{Rem} }
\newcommand{\rem}[1]{\begin{remark}{#1}}
\newcommand{\erem}{\end{remark}}

\newcommand{\eq}[1]%
        { \ifShowLabels \TeXrefEq{E:#1} \fi
           \begin{equation} \label{E:#1} }
\newcommand{\eeq}{\end{equation}}

\newcommand{\prf}{ \begin{proof} }
\newcommand{\eprf}{ \end{proof} }
\newcommand{\Label}[1]{\label{#1}  \ifShowLabels \TeXref{#1} \fi }

\ShowLabelsfalse

\newcommand{\n}{\nabla}
\renewcommand{\b}{\bullet}
\newcommand{\p}{\text{\( \partial\)}}
\newcommand{\F}{\calF}
\renewcommand{\tilF}{\widetilde{\calF}}
\newcommand{\T}{\calT}

\newcommand{\w}{\text{\( \ome\)}}

\newcommand{\we}{\text{\( \Ome^\b(\tilM,\tilF)\)}}

\renewcommand{\r}{\text{\( r^*_{\sqrt{T}}\)}}

\newcommand{\N}{{\tilN}}

\newcommand{\nb}{\nabla_T}
\newcommand{\lap}{\Del}

\renewcommand{\H}{H}
\newcommand{\E}{\calE}
\newcommand{\Ff}{\widetilde{\F}}

\renewcommand{\uc}{\tilU_c^-}

\newcommand{\ucp}{\tilU_{c'}^-}

\newcommand{\nce}{{\tilN_{C}}}

\newcommand{\Cech}{\v{C}ech\ }
\newcommand{\cS}{\mathcal{S} }

\renewcommand{\oM}{M}

\begin{document}

\title{Kirwan-Novikov inequalities on a manifold with boundary}
\author{Maxim Braverman}
\author{Valentin Silantyev}
\address{Department of Mathematics\\
        Northeastern University   \\
        Boston, MA 02115 \\
        USA
         }
\email{maxim@neu.edu} \email{v.silantyev@neu.edu}
\thanks{Research was partially supported by the NSF grant DMS-0204421}

\begin{abstract}
We extend the Novikov Morse-type inequalities for closed 1-forms in 2 directions. First, we
consider manifolds with boundary. Second, we allow a very degenerate structure of the critical
set of the form, assuming only that the form is non-degenerated in the sense of Kirwan. In
particular, we obtain a generalization of a result of Floer about the usual Morse inequalities
on a manifold with boundary. We also obtain an equivariant version of our inequalities.

Our proof is based on an application of the Witten deformation technique. The main novelty here
is that we consider the neighborhood of the critical set as a manifold with a cylindrical end.
This leads to a considerable simplification of the local analysis. In particular, we obtain a
new analytic proof of the Morse-Bott inequalities on a closed manifold.
\end{abstract} \maketitle

\tableofcontents
\setcounter{section}{-1}

\section{Introduction}\Label{S:introd}

\subsection{The main results}\Label{SS:results}
In \cite{BrFar1,BrFar2}, Michael Farber and the first author extended the Novikov Morse-type
inequalities for a closed 1-form  $\ome$ on a closed manifold to the case when  the critical set
of $\ome$ consists of a disjoint union of smooth submanifolds non-degenerate in the sense of
Bott. The Novikov inequalities were also considerably strengthen by means of twisting by an
arbitrary flat vector bundle $\calF$.

In this paper we further extend the results of \cite{BrFar1,BrFar2} by considering a manifold
$\oM$ with boundary and by allowing much more degenerate critical sets of $\ome$. For our
inequalities to be valid  we have to assume that near the boundary the 1-form  is exact and can
be written as
\footnote{In fact, our results are valid under a slightly weaker assumption. See \refs{prelim} for
the exact formulation.}
\eq{w=dtf}
 \ome \ =  \ d\, \big(\, f(x)t^m/m\,\big), \qquad m>0,
\end{equation}
where $x$ is a coordinate on the boundary $\Gam= \p\oM$ of $\oM$,  $t>0$ is an additional
coordinate such that $\Gam= \{t=1\}$, and $f$ is a smooth function on $\Gam$ for which zero is a
regular value. Set
\[
    U^- \ :=\ \big\{\, x\in\Gam:\, f(x)<0\, \big\}.
\]
If the form $\ome$ is exact than the role of the Betti numbers in our Morse-type inequalities is
played by the dimensions of the relative cohomology $H^\b(\oM,U^-;\calF)$ of the topological
pair $(\oM,U^-)$ with coefficients in a flat vector bundle $\calF$. In general, we define the
{\em relative Novikov numbers  of the topological pair $(\oM,U^-)$ with coefficients in a flat
vector bundle $\calF$}, which enter our inequalities, cf. \refd{novnum}. Our construction of the
relative Novikov numbers is similar to the Pazhitnov's construction of the usual Novikov numbers
as generic values of the dimensions of the cohomology of a 1-parameter family of flat
connections \cite{Pazhitnov87} (see also \cite{BrFar2}).

When applied to an exact 1-form $\ome= dh$ our inequalities extend the Kirwan-Morse inequalities
to manifolds with boundary. For the case when the critical points of $\ome$ are isolated and
$\calF= \CC$ is a trivial line bundle, similar inequalities were obtained by Floer
\cite{Floer89}.

The condition \refe{w=dtf} is equivalent to the assumption that $\ome$ can be extended to a
homogeneous 1-form on the manifold
\eq{tilM=}
    \tilM \ = \ \oM \cup \Big(\, \Gam\times[1,\infty)\,\Big),
\end{equation}
which does not have zeros on $\Gam\times[1,\infty)$, cf. \refs{cylend}. We refer to $\tilM$ as a
{\em manifold with a cylindrical end}. Thus we also obtain Morse-type inequalities for
homogeneous 1-forms on a manifold with a cylindrical end, which generalize the Morse
inequalities for generating functions quadratic at infinity, used in the theory of Lagrangian
intersections \cite{Viterbo92,EliashGromov98}.

We also obtain an equivariant version of our inequalities, which extends the results of
\cite{BrFar3,BrFar4} to manifolds with boundary and to 1-forms whose critical sets are not
manifolds.  In particular we obtain equivariant Morse inequalities on a manifold with boundary.
It is also important that we consider Morse-type inequalities twisted by a flat vector bundle
$\calF$. In many examples, this gives much stronger inequalities than the usual Morse-type
inequalities with coefficients in $\CC$, cf. Example~1.7 of \cite{BrFar2}. This is even more
important for equivariant inequalities, cf. \cite{BrFar3} and also the discussion in
\refs{equiv}.

\subsection{The method of the proof}\Label{SS:methodproof}
Our proof is based on an application of a version of the Witten deformation technique applied to
the extension of all the structures to the manifold \refe{tilM=}. Though a large part of our
proof is quite standard, there are two new ideas involved. First, the structure of our
Witten-type deformation near the critical set of $\ome$ is new. To explain the novelty, let us
temporary assume that $\ome$ is non-degenerate in the sense of Bott. Recall that the Witten
technique is based on considering a one-parameter deformation $\Del_T$ $(T\in \RR)$ of the
Laplacian and comparing the spectrum of $\Del_T$ for large $T$ with the spectrum of a similar
operator $\Del^N_T$ on the normal bundle to the critical set of $\ome$. It is crucial for the
method to work that the eigenfunctions of both operators $\Del_T$ and $\Del_T^N$ concentrate
near the critical set for large $T$. It is well known, however, that if the critical points of
$\ome$ are not isolated and if we endow the bundle $N$ with a natural ``bundle-like" Riemannian
metric then the naive generalization of the Witten construction does not work, because the
eigenfunctions of $\Del_T^N$ do not concentrate to the critical set as $T\to\infty$, cf.
\cite{Bis1}. This problem was solved by Bismut \cite{Bis1}, who proposed a smart two-parameter
family of deformations of the Laplacian (the Bismut deformation). It is not clear, however, how
to define this two-parameter family when the critical set of $\ome$ is not a manifold.

In this paper we apply a different idea.  We consider the bundle $N$ as a manifold with a
cylindrical end and introduce a Riemannian metric on $N$ which is conical on the cylindrical
end, cf. Sections~\ref{S:witten} and \ref{S:witN}. One of the advantages of this approach is
that now the analysis of the spectrum of the operator $\Del_T^N$ is exactly the same as the
analysis of the spectrum of the operator $\Del_T^{\tilM}$ on the manifold $\tilM$. In
particular, the eigenfunctions of $\Del_T^N$ concentrate near the critical set and the kernel of
$\Del_T^N$ can be calculated in cohomological terms, cf. Sections~\ref{S:witten} and
\ref{S:prkerH=}.

The second new element of our proof is the study of the Witten Laplacian on a manifold with a
cylindrical end (which we need for the study of both operators $\Del_T^N$ and $\Del_T^{\tilM}$).
In \refs{prkerH=}, we calculate the dimension of the kernels of the Witten Laplacians on a
manifold with a cylindrical end and show that they are equal to the relative Novikov numbers  of
the topological pair $(\oM,U^-)$ with coefficients in a flat vector bundle $\calF$. This is
probably the most non-trivial part of our paper.

\subsection{The structure of the paper}\Label{SS:strpaper}
In \refs{prelim} we formulate our main result and discuss some of its applications and
implications.

In \refs{equiv} we present an equivariant version of our inequalities.

In \refs{cylend} we reformulate our main theorem in terms of a manifold with a cylindrical end
and show the equivalence of the two formulations.

In \refs{witten} we define the Witten-type deformation of the Laplacian on a manifold with a
cylindrical end and show that it has a discreet spectrum.

In \refs{prkerH=} we show that the dimension of the kernel of the Witten Laplacian is given by
the relative Novikov numbers.

In \refs{witN} we introduce a structure of a manifold with a cylindrical end on a neighborhood
of the critical set of $\ome$ and describe the structure of our Witten-type deformation near the
critical set.

In \refs{proof} we present a proof of the Kirwan-Novikov inequalities, based on a result about a
comparison between the spectrum of the Witten Laplacians on $\tilM$ and on $N$. This result is
proven in \refs{model}.

\subsection*{Acknowledgment}
Some of the ideas used in this paper were developed several year ago in a joint unfinished
project of M.~Farber, M.~Shubin and the first author. We are very grateful to M.~Farber and
M.~Shubin for helping to develop these ideas and for allowing us to use them in this paper.




\section{Preliminaries and the main result} \Label{S:prelim}

\subsection{}\Label{SS:setting}
Let $\oM$ be a compact manifold with boundary $\Gam= \p\oM$.  Note that we don't exclude the
case when $\Gam$ is empty. If the boundary is not empty we will identify its tubular
neighborhood $\calU$ with the product $\Gam\times(0,1]$ and we will identify points of $\calU$
with pairs $(x,t), \ x\in\Gam, \, t\in (0,1]$.

Let $\ome$ be a closed 1-form on $M$ such that the restriction of $\ome$ to $\calU$ is exact. In
other words we assume that there exists a smooth function $h:\Gam\times(0,1]\to \RR$ such that
\eq{omeU=dh}
    \ome(x,t) \ = \ dh(x,t) \ = \ d_\Gam h(x,t)+\frac{\p h}{\p t}(x,t)\,dt,
    \qquad x\in\Gam, \ t\in (0,1].
\end{equation}
Here $d_\Gam$ is the de Rham differential on $\Gam$.

\subsection{Assumptions on $\ome$ at the boundary}\Label{SS:assbound}
We assume that
\begin{itemize}
    \item[\textbf{(B1)}] $d_\Gam{h}(x,1)+ \frac{\p h}{\p t}(x,1)\,dt \not= 0$ for all
    $x\in \Gam$.
    This implies that $\ome$ does not have zeros at a neighborhood of $\p{M}$. Without
    loss of generality we can and we will assume that {\em $\ome$ does not have zeros at
    $\Gam\times(0,1]$}.

    \item[\textbf{(B2)}] $0$ is a regular value of the function $\frac{\p h}{\p t}(x,1)$.

    \item[\textbf{(B3)}] If $\frac{\p h}{\p t}(x,1)= 0$ then
    \eq{dGnotdGt}
        d_\Gam h(x,1) \ \not=  \ -\lam\cdot d_\Gam\, \Big(\, \frac{\p h}{\p t}(x,1)\, \Big),
        \qquad\text{for all}\quad \lam>0.
    \end{equation}
\end{itemize}

The assumptions (B1) and (B2) mean that $\ome$ is in a generic position near $\p{\oM}$. They can
be always achieved by a small perturbation of $h$. Note also that (B2) implies, in particular,
that the open set
\eq{U-}
    U^- \ = \ \big\{\, x\in \Gam:\, \frac{\p h}{\p t}(x,1)<0\, \big\}\times\{1\} \subset \p\oM
\end{equation}
has a smooth boundary.

The assumption (B3) is more restrictive. It is essentially equivalent to the fact that the
closure $\overline{U^-}$ of $U^-$ is the {\em exit set} in the sense of Conley
\cite{Conley78,Floer89} for the flow of the vector field $-\n{}h$.  In \refs{cylend} we also
explain that this condition is necessary to extend $\ome$ to a ``homogeneous at infinity" form
on $\oM\cup \big(\Gam\times(1,\infty)\big)$. Also in \refss{exB3} we give an example which shows
that without this assumption our main theorem (\reft{main}) is not valid.

\subsection{Example}\Label{SS:exft2}
Suppose that on near $\Gam$ we have $h(x,t)= \frac1m{}f(x)t^m,\,  m>0$. Then condition (B3) is
automatically satisfied, while conditions (B1) and (B2) hold if and only if 0 is a regular value
of $f(x)$. This example is the most important for applications. In fact, in \refs{cylend} we
will show that if the conditions (B1)-(B3) are satisfied, then one can always extend $\ome$ and
$h$ to a bigger manifold $\oM'= \oM\cup(\Gam\times[1,A])$ without producing new zeros so that
near the boundary of $\oM'$ we have $h(x,t)= f(x)t^2/2$.

\subsection{Non-degeneracy assumptions on the zeros of $\ome$} \Label{SS:nondeg}
We assume that the set
\[
    \bfC \ =\ \big\{\,x\in M:\,\, \ome(x)=0\,\big\},
\]
called the {\em critical set} of $\ome$, belongs to $\oM\backslash\calU$ and satisfies the
following {\em non-degeneracy condition} of Kirwan \cite{Kirwan84}:

\defe{kirwan}
A closed 1-form $\ome$ is called {\em minimally degenerate} or {\em non-degenerate in the sense
of Kirwan} if it satisfies the following conditions
\begin{itemize}
\item[\textbf{(C1)}]
The critical set $\bfC$ is a finite union of disjoint closed subsets $C\subset \bfC$ called {\em
critical subsets of $\ome$}. In a neighborhood of each subset $C$ there exists a smooth function
$h_C$ such that $\ome=dh_C$ and $h_C(x)=0$ for all $x\in C$.

\item[\textbf{(C2)}] For every $C\subset \bfC$ there exists a locally closed connected submanifold
$\Sig_C$, called a {\em minimizing manifold}, containing $C$ such that $h_C|_{\Sig_C}\ge 0$ and
\eq{C2}
    C \ = \ \big\{\, x\in \Sig_C:\, h_C(x)=0\, \big\}.
\end{equation}
In other words, $C$ is the subset of $\Sig_C$ on  which $h_C$ takes its minimum values.

\item[\textbf{(C3)}]
At every point $x\in C$ the tangent space $T_x\Sig_C$ is the maximal among all subspaces of
$T_xM$ on which the Hessian $H_x(h_C)$ is positive semi-definite.
\end{itemize}
\edefe
The minimal degeneracy means that critical subsets can be as degenerate as a minimum of a
function, but not worse.

\rem{C3}
The condition (C3) can be reformulated as follows. Let $p:\nu(\Sig_C)\to \Sig_C$ be the normal
bundle to $\Sig_C$ in $M$. Fix a Euclidian metric on $\nu(\Sig_C)$ and let $|y|$ denote the norm
of a vector $y\in \nu(\Sig_C)$ with respect to this metric. Then by the generalized Morse lemma
\cite[Ch.~6]{Hirsch}  the condition (C3) is equivalent to the existence of a neighborhood
$\calW$ of the zero section in $\nu(\Sig_C)$ and of an embedding $i:\calW\to M$ such that
\eq{f(y)}
    (h_C\circ i)(y) \ = \ h_C(p(y)) \ - \ \frac{|y|^2}2
\end{equation}
and $h_C$ does not have critical points on $\calW\backslash C$.
\erem
\rem{nocrpt}
By the condition (C1), for every $C\in \bfC$ there exists a neighborhood $\calU$ of $C$ such
that $\ome$ does not have critical points in $\calU\backslash{}C$. Hence, it follows from
\refe{f(y)} that the restriction of $h_C$ to $\Sig_C$ also does not have critical points on
$\calU\cap\Sig_C$.
\erem
\subsection{The Morse counting polynomial} \Label{SS:calM}
The dimension of the fibers of $p:\nu(\Sig_C)\to \Sig_C$ is called the {\em index} of the
critical subset $C$ and is denoted  by $\ind(C)$.

Let $o(\Sig_C)$ be the {\em orientation bundle of $\nu(\Sig_C)$, considered as a flat line
bundle} over $\Sig_C$. We denote by $o(C)$ the restriction of $o(\Sig_C)$ to $C$. Let $\F$ be a
flat vector bundle over $M$. Consider the {\em twisted Poincar\'e polynomial} of the critical
subset $C$
\begin{equation}\Label{E:P-pol}
    \calP_{C,\F}(\lambda)\ =\
       \sum_{i=0}^{\dim M}\, \lambda^i\dim_{\CC} \check{H}^i(C,\F_{|_C}\otimes o(C)),
\end{equation}
where $\check{H}^i(C,\F_{|_C}\otimes o(C))$ denotes the \Cech cohomology of $C$ with
coefficients in the flat vector bundle $\F_{|_C}\otimes o(C)$. Define the following {\em Morse
counting polynomial}
\begin{equation}\Label{E:M-pol}
    \calM_{\omega,\F}(\lambda)\ =
        \ \sum_{C\in \bfC} \lambda^{\ind(C)} \calP_{C,\F}(\lambda),
\end{equation}
where the sum is taken over all critical subsets $C\subset \bfC$.

\subsection{The generalized Novikov numbers} \Label{SS:calN}
Assume that $\F$ is a complex flat finite dimensional vector bundle over $M$ and let
\[
        \nabla:\, \Omega^\b(M,\F)\ \to\ \Omega^{\b+1}(M,\F)
\]
denote the covariant differential on $\F$.

A  closed  1-form $\omega\in \Omega^1(M)$ on $M$ with real values determines a family of
connections on $\F$ (the {\em Novikov deformation}) parameterized by the real numbers $T\in \RR$
\begin{equation}\Label{E:nab-t}
    \nabla_T:\, \Omega^i(M,\F)\ \to\ \Omega^{i+1}(M,\F); \qquad
        \nb=\n+Te(\w).
\end{equation}
Here $e(\w)$ denotes the operator of exterior multiplication by \w. All the connections
$\nabla_T$ are flat, i.e., $\nabla_T^2=0$, if the form $\omega$ is closed. Hence, for any $T\in
\RR$ the pair $(\F,\n_T)$ is a flat vector bundle over $M$. Sometimes  we will denote this flat
bundle by $\F_T$ for short.

The flat bundle $\F_T$ admits the following alternative description. For $T\in\RR$, let
$\calE_{T\ome}$ denote the flat real line bundle over $M$ with the monodromy representation
$\rho_{T\ome}:\pi_1(M)\to \RR^\ast$ given by the formula
\begin{equation}\Label{E:novrep}
    \rho_{T\ome}(\gamma)=
        \exp(-T\int_\gamma\omega)\ \ \in \ \RR^\ast,\qquad \gamma\in
                            \pi_1(M).
\end{equation}
Then $\F_T$ is isomorphic to the tensor product $\F\otimes\calE_{T\ome}$.

Assume now that $\w$ is a non-degenerate in the sense of Kirwan closed 1-form and let $U^-$ be
as in \refe{U-}. Let $H^\b(\oM,U^-;\F_T)$ be the relative cohomology of the pair $(\oM,U^-)$
with coefficients in the flat bundle $\F_T$.


The dimension of the cohomology $H^\b(\oM,U^-;\F_T)$ is an integer valued function of $T\in\RR$.
This function has the following behavior: There exists a discrete subset $S\subset\RR$, such
that the dimension $\dim H^\b(\oM,U^-;\F_T)$ is {\it constant for $T\notin S$} (the
corresponding value of the dimension we will call the {\em background value}; the corresponding
value of $T$ is called {\em generic}) and for $T\in S$ the dimension of $H^\b(\oM,U^-;\F_T)$ is
{\em greater} than the background value. Cf., for example, \cite[Theorem 2.8]{Farber85}, where a
more precise information for the case of general elliptic complexes is given. The subset $S$
above will be called {\em the set of jump points}.

\defe{novnum} For each $i=0,1\nek n$, the background value of the
  dimension of the cohomology $H^i(\oM,U^-;\F_T)$ is called the $i$-th
  {\em generalized Novikov number} and is denoted by $\beta_i(\ome,\F)$.
\edefe
\rem{ome1-ome2}
Note that if the boundary of $M$ is empty then $\bet_i(\ome,\F)$ are the usual Novikov numbers
\cite{Novikov81,Novikov82,Pazhitnov87} with coefficients in $\calF$, cf. \cite{BrFar1,BrFar2}.
In this case $\bet_i(\ome,\F)$ depends only on the cohomology class of $\ome$. More generally,
let $\ome_1$ and $\ome_2$ be closed 1-form representing the same cohomology class in
$H^1(M,\CC)$ and satisfying the assumptions (B1)-(B3) of \refss{assbound}. Let $U^-_1,\ U^-_2$
be the sets defined by \refe{U-} using the forms $\ome_1$ and $\ome_2$ respectively. If the
topological pairs $(\oM,U^-_1)$ and $(\oM,U^-_2)$ are homotopically equivalent then
$\bet_i(\ome_1,\F)= \bet_i(\ome_2,\F)$.
\erem

We define the  {\em Novikov type polynomial} by the formula
\begin{equation}\Label{E:N-pol}
    \calN_{\ome,\F}(\lambda)\ =\ \sum_{i=0}^{\dim M}\, \lambda^i \beta_i(\ome,\F),
\end{equation}

The main result of this paper is the following Morse-type inequalities for a differential 1-form
on a manifold with boundary:
\th{main} Assume that $\w$ is a closed 1-form on $M$ which satisfies the assumptions
 (B1)-(B3) of \refss{assbound} and all whose zeros are non-degenerate in the sense of
 Kirwan, cf. \refss{nondeg}. Then there exists a polynomial
  $\calQ(\lambda)= q_0 +  q_1\lambda + q_2\lambda^2 + \dots$
  with non-negative integer coefficients $q_i\ge 0$, such that
  \begin{equation}\Label{E:main}
    \calM_{\omega, \F}(\lambda)\ -\ \calN_{\ome,\F}(\lambda)\ =\
            (1+\lambda)\calQ(\lambda).
  \end{equation}
\eth

The proof of the theorem occupies most of the remaining part of the paper.

\rem{BrFar}
Another version of Morse-type inequalities for differential 1-forms on manifolds with boundary
were suggested by M.~Farber and the first author in \cite{BrFar3}. The conditions on the
behaviour of $\ome$ near the boundary was, however, quite different in \cite{BrFar3}. In this
sense \reft{main} complements the results of \cite{BrFar3}. Note also that in \cite{BrFar3} the
zeros of the form were assumed to be non-degenerate in the sense of Bott, while in \reft{main}
we allow a much more complicated structure of the set of zeros.
\erem

\subsection{Corollaries and applications}\Label{SS:appl}
We will now discuss some special cases of \reft{main}. First, consider the case when the
boundary of $\oM$ is empty. Then $\bet_i(\ome,\calF)$ are the usual Novikov numbers with
coefficients in $\calF$ associated to the cohomology class of $\ome$, cf. \cite{BrFar2}.
\reft{main} in this case generalizes the classical Novikov inequalities in two directions:
first, we allow non-isolated zeros of $\ome$ (assuming they are non-degenerate in the sense of
Kirwan). And, second, we strengthen the inequalities by twisting them with the vector bundle
$\F$ (cf. example~1.7 in \cite{BrFar2}, which shows that the obtained inequalities are really
stronger than the classical inequalities with coefficients in $\CC$. We also refer to
\cite{BrAlesker} where the Kirwan inequalities with coefficients in a non-trivial bundle are
used to obtain a topological information about the symplectic reduction of a Hamiltonian
$T$-space). For the case when the zeros of $\ome$ are non-degenerate manifolds (i.e., when
$\ome$ satisfies Bott non-degeneracy conditions) these results were obtained by M.~Farber and
the first author in \cite{BrFar1,BrFar2}.

Let us return to the case when $\oM$ has boundary. Assume that the cohomology class of $\ome$
vanishes, i.e.,  there exists a  function $h:\oM\to \RR$ such that $\ome =dh$. Then
$\bet_i(\ome,\F)$ are just the dimensions of the cohomology of the pair $(\oM,U^-)$ with
coefficients in $\F$. Thus our result generalizes the classical Morse inequalities to the case
of a manifold with boundary.  For the special case when $h$ is a Morse function (i.e., has only
isolated non-degenerate critical points) and $\F= \CC$ is the trivial bundle, similar
inequalities were obtained by Floer \cite{Floer89}.

In \refs{equiv} we use \reft{main} to obtain the equivariant version of the Kirwan-Novikov
inequalities on a manifold with boundary. Note, that the use of cohomology with coefficients in
a flat vector bundle $\calF$ (rather than just complex values cohomology) is very important for
this application, cf. \refr{coefF} for details.

\subsection{The necessity of the condition (B3)}\Label{SS:exB3}
We finish this section with a very simple example, which shows that without the assumption (B3)
(cf. \refss{assbound}) \reft{main} is not valid.

Consider the cylinder $\oM= S^1\times [1,2]$. Then $\Gam = \p{}M= (S^1\times\{1\})\sqcup
(S^1\times\{2\})$. Set
\eq{hexB3}
    h(x,t) \ = \ t^2\sin x,
\end{equation}
where $x\in S^1, \ t\in [1,2]$. Let $\ome =dh$. One easily checks that all the conditions of
\reft{main} except (B3) hold. The condition (B3) fails at the points $(0,1)$ and
$(\pi,1)$.%
\footnote{Though \refe{hexB3} formally looks like the function in the Example~\ref{SS:exft2}, the
 situation here is quite different. This is because near $S^1\times\{1\}$ the parameter $t$ in
 \refe{hexB3} is not the same as in \refss{setting}.}

Clearly,
\[
    U^- \ =\ \big(\,(0,\pi)\times\{1\}\,\big)\sqcup \big(\,(\pi,2\pi)\times\{2\}\,\big).
\]
Let $\F=\CC$ be the trivial line bundle. Since the form $\ome$ is exact,
$\bet_i=\bet_i(\ome,\F)$ equals to the dimension of the relative cohomology $H^i(\oM,U^-)$. An
easy calculation shows that
\[
    \bet_0 \ = \ \bet_2 \ = \ 0, \ \ \bet_1 \ = \ 2.
\]
So $\calN_{\ome,\F}(\lam)=2\lam$ in this case. But $\ome$ does not have zeros on $M$. Hence, the
Morse counting polynomial is equal to 0 and the left hand side of \refe{main} equals $-2\lam$.
Thus \refe{main} can not hold.

Let us enlarge the manifold $\oM$ by setting $\oM'= S^1\times[-2,2]$. The form $\ome=dh$ extends
naturally to the bigger manifold. Now the condition (B3) is easily seen to be satisfied. We
still have $\calN_{\ome,\F}(\lam)= 2\lam$. But now $\ome$ vanishes on $S^1\times\{0\}\subset
S^1\times[-2,2]$. The Morse counting polynomial equals $2\lam$ and \refe{main} holds with
$\calQ=0$.

This example illustrates an important phenomenon to be discussed in the next section: the
condition (B3) allows a natural extension of $\ome$ from $\calU=\Gam\times(0,1)$ to the infinite
cylinder $\Gam\times(0,\infty)$ without producing new zeros of $\ome$.

\section{Equivariant Kirwan-Novikov inequalities on a manifold with boundary}\Label{S:equiv}

In this section we present a generalization of the equivariant Novikov inequalities obtained  by
M.~Farber and the first author in \cite{BrFar3,BrFar4} to manifolds with boundary and to 1-forms
which are not non-degenerate in the sense of Bott (but are non-degenerate in the sense of
Kirwan). The constructions and the proofs are very similar to \cite{BrFar4}, so we will be very
brief in this section and we will refer the reader to \cite{BrFar4} for details.

Throughout this section $G$ is a compact Lie group and $\oM$ is a compact $G$-manifold with
boundary $\Gam=\p\oM$.

\subsection{Basic 1-forms}\Label{SS:basicform}
Recall that a smooth 1-form $\ome$ on a $G$-manifold $M$ is called {\em basic} (cf.
\cite{AtBott84}) if it is $G$-invariant and its restriction on any orbit of the action of $G$
equals to zero.

Note that, if the group $G$ is finite, then $\ome$ is basic if and only if it is $G$-invariant,
i.e., if $g^*\ome=\ome$ is satisfied for every $g\in G$. Also, {\em if $M$ is connected and if
the set of fixed points of the action of $G$ on $M$ is not empty, then any closed $G$-invariant
1-form on $M$ is basic}, cf.  Lemma~3.4 of \cite{BrFar4}. Note also that {\em any exact
invariant form $\ome= df$ is basic}.

\subsection{The equivariant Novikov numbers}\Label{SS:equivnov}
Let $\calF\to \oM$ be a $G$-equivariant flat vector bundle over $\oM$, cf. Section~2 of
\cite{BrFar4}. Suppose $\ome$ is a closed basic 1-form on $M$ which satisfies assumptions
(B1)-(B3) of \refss{assbound}. Let $U^-\subset \Gam$ be as in \refe{U-}. We denote by
$H^\b_G(\oM,U^-;\calF)$ the equivariant cohomology of the pair $(\oM,U^-)$ with coefficients in
$\calF$.
\footnote{The equivariant cohomology $H_G^\b(M;\calF)$ with coefficients in an equivariant flat
bundle are defined, for example, in \S2 of \cite{BrFar4}. The relative equivariant cohomology
with coefficients is defined in exactly the same way.}

Every closed basic 1-form $\ome$ determines an equivariant flat line bundle $\E_\ome$ (with
$\ome$ being its connection form), cf. \refss{calN} (see also \S2.2 of \cite{BrFar4}). Using
these constructions we define the equivariant Novikov numbers as follows.

Given an equivariant flat bundle $\F$ over $M$, consider the one-parameter family $\F_T=
\F\otimes \E_{T\ome}$ of equivariant flat bundles, where $T\in \RR$, ({\em the Novikov
deformation}) and consider the twisted equivariant cohomology
\[
    H^i_G(\oM,U^-;\F\otimes\E_{T\ome}), \qquad\text{where}\quad t\in\RR,
\]
as a function of $T\in \RR$. The same arguments as in the proof of Lemma~1.5 of \cite{BrFar4}
show that, for each $i=0,1,\dots$, there exists  a finite subset $S\subset\RR$ such that the
dimension of the cohomology $H^i_G(\oM,U^-;\F\otimes\E_{T\ome})$ is constant for $T\notin S$ and
the dimension of the cohomology $H^i_G(\oM,U^-;\F\otimes\E_{T\ome})$ jumps up for $T\in S$.

The subset $S$, is called the {\em set of jump points}; the value of the dimension of
$H^i_G(\oM,U^-;\F\otimes\E_{T\ome})$ for $T\notin S$ is called {\em the background value of the
dimension of this family}.

\defe{equivnov}
The $i$-dimensional {\em equivariant Novikov number} $\bet_i^G(\ome,\F)$ is defined as the
background value of the dimension of the cohomology of the family
$H^i_G(\oM,U^-;\F\otimes\E_{T\ome})$.
\edefe

\rem{ome1-2}
If the boundary of $\oM$ is empty then the numbers $\bet_i^G(\ome,\F)$ coincide with the
equivariant Novikov numbers with coefficients in $\calF$ defined in \cite{BrFar4}. In this case
$\bet_i^G(\ome,\F)$ depends only on the cohomology class of $\ome$. More generally,
$\bet_i(\ome_1,\F)= \bet_i(\ome_2,\F)$ if $[\ome_1]= [\ome_2]\in H^1_G(M)$ and the topological
pairs $(\oM,U^-_1)$ and $(\oM,U^-_2)$ are homotopically equivalent (here $U_1^-,\, U_2^-\in
\Gam$ are the sets defined by \refe{U-} using the forms $\ome_1$ and $\ome_2$ respectively).
\erem

The formal power series
\[
        \calN_{\ome,\F}^G(\lam) \ =\
                \sum_{i=0}^\infty\,\lam^i\bet_i^G(\ome,\F)
\]
is  called {\em the equivariant Novikov series}.

\newcommand{\fo}{{\calF_{|_C}\otimes o(C)}}
\subsection{The equivariant Morse series}\Label{SS:equivMorse}
We assume that $\ome$ is non-degenerate in the sense of Kirwan and we use the notation
introduced in Subsections~\ref{SS:nondeg} and \ref{SS:calM}. Let $C\subset \bfC$ be a component
of the critical set of $\ome$. If the group $G$ is connected then $C$ is a $G$-invariant subset
of $M$. In general case we denote by
\[
    G_C \ := \ \big\{\, g\in G:\, g\cdot C\subset C\,\big\}
\]
the stabilizer of $C$ in $G$. Let $|G:G_C|$ denote the index of $G_C$ as a subgroup of $G$.
Since $G_C$ contains the connected component of the unity in $G$, this index is finite.

The compact Lie group $G_C$ acts on $C$ and the flat vector bundles $\F_{|_C}$ and $o(C)$ (cf.
\refss{calM}) are $G_C$-equivariant.  Let
 \(
        \check{H}_{G_C}^\ast(C,\fo)
 \)
denote the {\em equivariant \Cech cohomology} of the flat $G_C$-equivari\-ant vector bundle
$\fo$. Consider the {\em equivariant Poincar\'e series} of $C$
\[
        \calP^{G_C}_{C,\F}(\lam)\ =\
           \sum_{i=0}^{\infty}\, \lam^i\, \dim_{\CC}H_{G_C}^i(C,\fo)
\]
and define using it the following {\em equivariant Morse counting series}
\[
        \calM_{\w,\F}^G(\lam)\ = \ \sum_C\,
             \lam^{\ind(C)}\, |G:G_C|^{-1}\,\calP_{C,\F}^{G_C}(\lam),
\]
where the sum is taken over all components $C$ of $\bfC$.

The main result of this section is the following:
\th{equiv}
Suppose that $G$ is a compact  Lie group and $\F$ is a flat $G$-equivariant vector bundle over a
compact $G$-manifold $M$ with boundary. Let $\w$ be a closed basic 1-form on $M$ which satisfies
the assumptions (B1)-(B3) of \refss{assbound} and all whose zeros are non-degenerate in the
sense of Kirwan, cf. \refss{nondeg}. Then there exists a formal power series $\calQ(\lam)$ with
non-negative integer coefficients, such that
\[
    \calM_{\w,\F}^G(\lam) \ -\ \calN_{\xi,\F}^G(\lam)\ =\
                (1+\lam)\,\calQ(\lam).
\]
\eth
The theorem follows from \reft{main} by exactly the same arguments as were used in \cite{BrFar4}
to derive the equivariant Novikov inequalities for a non-degenerate in the sense of Bott one
form on a closed manifold from the Novikov-Bott inequalities of \cite{BrFar2}. We leave the
details to the interested reader.

\rem{coefF}
Let us remark that the use of the cohomology twisted by a flat vector bundle $\calF$
considerably strengthen the inequalities of \reft{equiv}. For example, assume that the group $G$
is finite and the form $\ome$ is exact. An application of the equivariant Morse inequalities of
Atiyah and Bott \cite{AtBott82,Bott2} leads to estimates which are often weaker than the
standard Morse inequalities (ignoring the group action). However, the application of the
inequalities twisted by an equivariant flat bundle leads to a much better estimates, cf.
\cite{BrFar3}. This was used in \cite{BrAlesker}, where an application of \reft{equiv} to a
finite group but with a non-trivial equivariant vector bundle provided a new topological
information about the cohomology of the symplectic reduction of a Hamiltonian $T$-manifold.
\erem


\section{Reformulation in terms of a manifold with a cylindrical end}\Label{S:cylend}

In this section we reformulate \reft{main} as Morse-type inequalities for differential forms on
a manifold with a cylindrical end. This new formulation is more suitable for our analytical
proof, but also has an independent interest. In particular, it generalizes the Morse
inequalities for generating functions quadratic at infinity, which are used in the theory of
Lagrangian intersections \cite{Viterbo92,EliashGromov98}.


\subsection{A manifold with a cylindrical end}\Label{SS:cylindrical}
By a  {\em manifold with a cylindrical end} we  understand a smooth (non-compact) manifold
$\tilM$ without boundary which has the form $\Gam\times (1,\infty)$ near infinity. More
precisely, we suppose that there exists a compact submanifold $\oM\subset \tilM$ with smooth
boundary $\Gam$ such that
\begin{equation}\Label{E:cyl}
    \tilM \ = \ \oM \cup \big(\Gam\times [1,\infty)\big)
\end{equation}
with $\p\oM$ being identified with $\Gam\times\{1\}$.

The submanifold $\T= \Gam\times (1,\infty)$ is called the {\em cylindrical end} of $M$.

Of course, $\oM$ is not unique. For example, it may be replaced by $\oM_{t_0}= \oM\cup
(\Gam\times (1,t_0])$. Similarly, the decomposition $\T$ into a direct product is not unique.
But we will fix $\oM$ and the decomposition $\T= \Gam\times (1,\infty)$ for the sake of
simplicity of the notation.

\subsection{Homogeneous 1-forms}\Label{SS:homog}
Let $\tilM=\oM\cup (\Gam\times [1,\infty))$ be a manifold with a cylindrical end.  Every point
of $\calT=\Gam\times (1,\infty)$ may be identified with a pair $(x,t)$ where $x\in \Gam,\ t>1$.

Let $\tau_s:\T\to \T$ $(s\ge 1)$ denote the multiplication by $s:\,\tau_s(x,t)=(x,s\cdot t)$.

A differential 1-form $\alp$ on $M$ is called {\em homogeneous of degree $m$ at infinity} if
there exists $t_0> 1$ such that
\begin{equation}\Label{E:hom}
    \tau_s^*(\alp_{|_{\Gam\times(t_0,\infty)}}) \ = \ s^m \alp_{|_{\Gam\times(t_0,\infty)}},
    \qquad \text{for every} \quad s>1.
\end{equation}

If $\ome$ is a homogeneous 1-form of degree $m$, then its restriction to
${\Gam\times(t_0,\infty)}$ may be written as
\[
    \alp_{|_{\Gam\times(t_0,\infty)}} \ = \ t^m a \ + \ f(x)t^{m-1}dt,
\]
where $a$ is a 1-form on $\Gam$. If, in addition,  $\alp$ is closed, $d\alp= 0$, then $a=
\frac1md{f}$ and
\begin{equation} \Label{E:w=df}
    \alp_{|_{\Gam\times(t_0,\infty)}} \ = \ d\,\Big(\,\frac1m f(x)t^m\,\big).
\end{equation}
In particular, {\em the restriction of a closed homogeneous at infinity 1-form to the
cylindrical part $\calT$ is always exact}.

\subsection{Non-degeneracy at infinity}\Label{SS:nondeginf}
Suppose that $\ome$ is a closed 1-form on $\tilM$ which is homogeneous at infinity. We say that
$\ome$ is {\em non-degenerate at infinity} if the restriction $\ome(x,t)$ of $\ome$ to the
cylindrical part $\calT$ does not have zeros for $t\gg1$.

Let $\ome$ be represented as in \refe{w=df}. Then $\ome$ is non-degenerate at infinity if and
only if 0 is a regular value of the function $f:\Gam\to \RR$.

\subsection{Example: a vector bundle over a compact manifold}\Label{SS:ex-vb}
One of the most important examples of a manifold with a cylindrical end is a vector bundle over
a compact manifold.

Let $C$ be a smooth compact manifold without boundary and let $E$ be a finite dimensional vector
bundle over $C$. Fix a Hermitian metric $h^E$ on $E$. For a vector $y\in E$ we denote by $|y|$
its norm with respect to the metric $h^E$.

Fix $r>0$ and define
\[
    B_r \ = \ \big\{\,y\in \E:\, |y|\le r\,\big\}; \qquad  S_r \ =\ \big\{\,y\in \E:\, |y|=r\,\big\}
\]
Then $B_r$ is a compact manifold, the boundary of $B_r$ equals $S_r$, and $E\backslash B_r$ is
diffeomorphic to the product $S_r\times (1,\infty)$. Hence, $E$ has a natural structure of a
manifold with a cylindrical end.

Assume now that the bundle $E$ splits into orthogonal direct sum $E=E^+\oplus{}E^-$. Consider
the function
\[
    h(y) \ = \ \frac{|y^+|^2}2- \frac{|y^-|^2}2, \qquad (y^+,y^-)\in E=E^+\oplus E^-.
\]
Then the form $\ome= dh$ is homogenous of degree 2 on $E$ and is non-degenerate at infinity.

More generally, if $\tilV$ is a manifold with a cylindrical end and $\tilN\to \tilV$ is a vector
bundle over $\tilV$, then one can introduce a structure of a manifold with a cylindrical end on
the total space of $\tilN$. A version of this example will allow us to consider the normal
bundle to a neighborhood of the critical set $C$ in the minimizing manifold $\Sig_C$ (cf.
\refd{kirwan}) as a manifold with a cylindrical end. We refer to \refss{tilN} for details.

\subsection{The Novikov numbers}\Label{SS:novnumb2}
Let $\ome$ be a closed 1-form on $\tilM$ which is homogeneous and non-degenerate at infinity.
Let $f(x)$ be as in \refe{w=df}.

Let $\tilF$ be a flat vector bundle over $\tilM$. As in \refss{calN} we construct a family
$\tilF_T \ (T\in \RR)$ of flat vector bundles over $\tilM$ using the 1-form $\ome$.

For every $c>0$ define
\eq{Uc}
    \uc = \big\{\, (x,t)\in \Gam\times(1,\infty):\, f(x)t^{m}/m<-c\, \big\}
\end{equation}
and consider the pair  $(\tilM,\uc)$.   It is clear, that if we change $c$ then this pair will
be replaced by a homotopic one. Thus the cohomology $H^i(\tilM,\uc;\tilF_T)$ of this pair with
coefficients in $\tilF_T$ is independent of the choice of $c$. We define the Novikov number
$\bet_i(\ome,\tilF)$ to be the {\em background value} of the dimension of this cohomology, cf.
\refss{calN}, and we define the Novikov polynomial $\calN_{\ome,\tilF}(\lam)$ by \refe{N-pol}.

If all the zeros of $\ome$ are non-degenerate in the sense of Kirwan we define the Morse
counting polynomial $\calM_{\ome,\tilF}(\lam)$ as in \refss{calM}.

\th{cylinder}
Let $\tilM$ be a manifold with cylindrical end and let $\tilF$ be a flat vector bundle over
$\tilM$. Let $\ome$ be a closed 1-form on $\tilM$ which is homogeneous and non-degenerate at
infinity. Assume also that all the zeros of $\ome$ are non-degenerate in the sense of Kirwan.
Then there exists a polynomial $\calQ(\lambda)$ with non-negative integer coefficients such that
the Morse-type inequalities \refe{main} hold.
\eth

The following proposition shows that this theorem is equivalent to \reft{main}.

\prop{cyl=bound}
\reft{cylinder} is equivalent to \reft{main}.  Moreover, \reft{cylinder} (and, hence,
\reft{main}) follows from its particular case when  $\ome$ is homogeneous of degree 2.
\eprop
Thus, to prove Theorems~\ref{T:main} and \ref{T:cylinder}, it is enough to prove \reft{cylinder}
for the case when  $\ome$ is homogeneous of degree 2. This will be done in the subsequent
sections. The rest of this section is occupied with the proof of \refp{cyl=bound}.

\subsection{\reft{main} implies \reft{cylinder}}\Label{SS:cyl<main}
Let $\ome$ be homogeneous, non-degenerate at infinity closed 1-form on a manifold $\tilM$ with
cylindrical end.

Fix a decomposition \refe{cyl} of $\tilM$ and let $t_0\ge 0$ be as in \refe{hom}. Set
\[
    \oM' \ = \ \oM\sqcup\big(\Gam\times(1,t_0+2]\big).
\]

Introducing a new coordinate $t'=t-t_0-1$ we identify a neighborhood $\calU$ of $\p\oM'$ with
the product $\Gam\times(0,1]$. From \refe{w=df} we conclude that on $\calU$ the form
$\ome=dh(x,t')$ where
\eq{ht'}
    h(x,t') \ = \ \frac1m f(x)(t'+t_0+1)^m.
\end{equation}
As in \refss{assbound} set
\[
    U^- \ = \ \big\{\, (x,t'):\, \frac{\p h}{\p t'}(x,1)<0,\, t'=1\, \big\}\ = \
        \big\{\, (x,t'):\, f(x)<0,\, t'=1\, \big\}.
\]

Recall that $\uc$ is defined in \refe{Uc}. Clearly, the pair $(\tilM,\uc)$ is homotopic to the
pair $(\oM',U^-)$. Thus the relative cohomology $H^\b(\tilM,\uc;\tilF_T)$ is naturally
isomorphic to $H^\b(\oM,U^-;\tilF_T|_{\oM'})$. \reft{cylinder} follows now from application of
\reft{main} to the form $\ome|_{\oM'}$ .\hfill$\square$

To show the implication in the other direction we need the following
\lem{tilh}
Let $h:\Gam\times(0,1]\to \RR$ be a smooth function satisfying the conditions (B1)-(B3) of
\refss{assbound}. According to the condition (B1) of \refss{assbound}, $dh\not=0$ on
$\Gam\times(0,1]$, cf. \refss{assbound}. Set
\[
    f(x) \ = \ \frac{\p h}{\p t}(x,1).
\]
Then there exist $\eps\in (0,1/3), \ m>0$ and a function $\tilh:\Gam\times(0,\infty)\to \RR$
such that
\begin{align}
    \tilh(x,t) \ = \ h(x,t) \qquad &\text{for}\quad t<1-2\eps, \ x\in\Gam; \Label{E:tilh=h}\\
    \tilh(x,t) \ = \ m\, f(x)\,\frac{t^2}2 \qquad &\text{for}\quad t>1-\eps, \ x\in \Gam;\Label{E:tilh=f}\\
    d\tilh(x,t) \ \not=\ 0 \qquad &\text{for all}\quad t>0, \ x\in \Gam. \Label{E:no0}
\end{align}
\elem
\prf
From the conditions (B3) of \refss{assbound} we see that there exist $\eps\in (0,1/3)$ such that
for all $t\in (1-3\eps,1]$ we have
\begin{align}
    d_\Gam h(x,t) \ &\not= \ -\lam d_\Gam f(x) \quad
    \text{for all} \ \lam>0 \ \ \text{and all} \ x\in\Gam \ \text{such that} \ |f(x)|< 2\eps. \Label{E:f<2eps} \\
    \frac{\p h(x,t)}{\p t} \ &\not= \ -\, f(x) \quad
    \text{for all} \ x\in\Gam \ \text{such that} \ |f(x)|>\eps. \Label{E:f>eps}
\end{align}

Set
\eq{t0}
    m \ = \ \frac{18}{\eps}\, \sup\big\{\, |h(x,t)|:\, x\in \Gam,\, 0<t\le 1\, \big\}.
\end{equation}

Let $\chi:\RR\to [0,1]$ be a smooth non-decreasing function such that $\chi(t)=0$ for $t\le
1-2\eps$ and $\chi(t)=1$ for $t\ge 1-\eps$. Set
\eq{tilh}
    \tilh(x,t) \ = \ (1-\chi(t))\, h(x,t)\ + \ \chi(t)\, m\,f(x)\,\frac{t^2}2.
\end{equation}

Clearly, $\tilh$ satisfies \refe{tilh=h} and \refe{tilh=f}. We claim that it also satisfies
\refe{no0}. Indeed,
\begin{itemize}
\item
If $0<t<1-2\eps$ then $\tilh(x,t)= h(x,t)$ and \refe{no0} is exactly our assumption on $h$.
\item
If $1-2\eps\le t\le 1$ and $|f(x)|<2\eps$, then, by \refe{f<2eps},
\[
    d_\Gam \tilh(x,t) \ = \
    (1-\chi(t))\, d_\Gam h(x,t)\ + \ \chi(t)\, \frac{m t^2}2 d_\Gam\,f(x) \ \not= \ 0.
\]
\item
If $1-2\eps\le t\le 1$ and $|f(x)|>\eps$, then
\[
    \frac{\p \tilh(x,t)}{\p t} \ = \ \chi'(t)\,\Big(\, m\,f(x)\,\frac{t^2}2-h(x,t)\,\Big) \ + \
    (1-\chi(t))\, \frac{\p h(x,t)}{\p t} \ + \ \chi(t)\, m\,f(x)\, t.
\]
By \refe{f>eps} and \refe{t0} all the terms in the right hand side of this equality have the
same sign as $f$. Hence, $\p\tilh/\p{t}\not=0$.
\item
If $t\ge 1-\eps$ then $h(x,t)=mf(x)t^2/2$. By assumption (B2) of \refss{assbound}, zero is a
regular value of $f(x)=\frac{\p{}h}{\p{t}}(x,1)$. Hence,
\[
    dh(x,t) \ = \ m\,f(x)\,t\, dt \ + \ \frac{m t^2}2\, d_\Gam f(x) \ \not= \ 0.
\]
\end{itemize}
The proof of the lemma is complete.
\eprf

\subsection{\reft{cylinder} implies \reft{main}}\Label{SS:cyl>main}
We are ready now to show that \reft{main} follows from the special case of \reft{cylinder} when
$\ome$ is homogeneous of degree 2 at infinity. This will finish the proof of \refp{cyl=bound}.

Let $\oM, \F$ and $\ome$ be as in \reft{main}. Consider the manifold
\[
    \tilM \ = \ \oM\, \cup\, \big(\,\Gam\times[1,\infty)\,\big)
\]
with a cylindrical end.

Fix a smooth diffeomorphism $s:(0,\infty)\to (0,1)$ such that $s(t)=t$ for $t\in (0,1/3]$.
Define a map  $\os:\tilM\to M$ such that $\os(m)=m$ for $m\in
\tilM\backslash\big(\Gam\times[1/3,\infty)\big)$ and $\os(x,t)= (x,s(t))$ for $x\in\Gam,\,
t\in[1/3,\infty)$.

Let $\tilF= \os^*\F$ be the pull-back of the bundle $\F$ to $\tilM$. We endow $\tilF$ with the
pull-back connection $\os^*\n$.

Let $h(x,t)$ be as in \refss{setting}. Assume, in addition, that $dh\not=0$ on
$\Gam\times(0,1]$, cf. \refss{assbound}. Let $\tilh$ be as in \refl{tilh}. Let $\tilome\in
\Ome^1(\tilM)$ be the 1-form whose restriction to
$\tilM\backslash\big(\Gam\times(1-2\eps,\infty)\big)$ is equal to $\ome$ and whose restriction
to $\Gam\times(1-2\eps,\infty)$ equals $d\tilh(x,t)$. From \refe{tilh=f} and \refe{no0} we
conclude that $\tilome$ is homogeneous of degree 2 and non-degenerate at infinity. Also, by
\refe{tilh=h}, the forms $\tilome$ and $\ome$ coincide in a neighborhood of their common set of
zeros. Hence their Morse counting polynomials coincide, $\calM_{\tilome,\tilF}(\lam)=
\calM_{\ome,\F}(\lam)$. The same arguments as in \refss{cyl<main} shows that
$\calN_{\tilome,\tilF}(\lam)= \calN_{\ome,\F}(\lam)$. Thus the Morse-type inequalities
\refe{main} for $\ome$ follow from application of \reft{cylinder} to the form $\tilome$. The
proof of \refp{cyl=bound} is complete. \hfill$\square$


\section{The Witten deformation of the Laplacian} \Label{S:witten}

Let $\tilM= \oM \cup \big(\Gam\times [1,\infty)\big)$ be a manifold with a cylindrical end and
let $\ome$ be a closed 1-form on $\tilM$ which is homogeneous of degree 2 at infinity, cf.
\refss{homog}.
Let $\tilF$ be a flat vector bundle over $\tilM$ endowed with a flat connection $\n$.

In this section we introduce a ``conical at infinity" metric on $\tilM$ and use it to construct
the Laplacian $\lap_T$ associated to $\nb$. This Laplacian plays a crucial role in our proof of
\reft{main}. In \refs{prkerH=}, we prove that the dimension of the kernel of the restriction of
$\lap_T$ to the space of $p$-forms is equal to the generalized Novikov number
$\bet_p(\xi,\tilF)$. In \refs{witN}, we estimate the dimension of the $\Ker\lap_T$ by comparing
it with the kernel of a certain Laplacian on a neighborhood of $\bfC$. That will prove
\reft{cylinder} (and, hence, in view of \refp{cyl=bound}, \reft{main}).

\subsection{Metrics on the cylindrical end} \Label{SS:metT}
The tangent bundle $T\T$ to the cylindrical end $\T=\Gam\times(1,\infty)$ splits into the
Whitney sum $T\Gam\oplus\RR$ of the tangent bundle $T\Gam$ to $\Gam$ and the tangent bundle to
the ray $(1,\infty)$ (which we identify with $\RR$).

Choose a Riemannian metric $g^\Gam$ on $\Gam$ and let $g^{\RR}$ denote the standard
(transitionally invariant) metric on $\RR$. Let $g^{\T}$ denote the metric on $\T$ given by the
formula
\begin{equation}\Label{E:gT}
    g^{\T} \ = \ t^2g^\Gam \ \oplus \ g^{\RR}.
\end{equation}
This metric is homogeneous of degree 2 in the following sense: Let $\tau_s:\T\to \T$ $(s\ge 1)$
denote the multiplication by $s:\,\tau_s(x,t)=(x,s\cdot t)$. Then
\eq{gThomog}
    \tau_s^*\, g^{\T} \ = \ s^2\, g^{\T}.
\end{equation}
We will refer to $g^{\T}$ as {\em conical metric} on $\T$.

We also fix a Hermitian metric $h^{\tilF|_\T}$ on the restriction $\tilF|_\T$ of the bundle
$\tilF$ to $\T$ which is flat along the ray $\{(x,t): \, t>1\}\subset \T$ for any $x\in \Gam$.

\subsection{Metrics on the manifold} \Label{SS:metM}
Let $g^\tilM$ be a Riemannian metric on $\tilM$ whose restriction on the cylindrical end $\calT$
is equal to the conical metric $g^\T$ (cf. \refss{metT}).

We also fix a Hermitian metric $h^{\tilF}$ on $\tilF$ whose restriction to $\T$ is equal to
$h^{\tilF|_\T}$.

\subsection{The deformation of the Laplacian} \Label{SS:witM}
Let $\nb=\n+Te(\w)$ be the Novikov deformation of the covariant derivative, cf. \refss{calN}.
Denote by $\nb^*$ the formal adjoint of $\nb$ with respect to the metrics $g^{\tilM},h^{\tilF}$.
\defe{laplM}
The {\it Witten Laplacian} is the operator
\[
    \lap_T\ := \ \frac1T(\nb\nb^*+\nb^*\nb):\, \Ome^\b(\tilM,\tilF)\, \to\,
    \Ome^\b(\tilM,\tilF).
\]
\edefe
It is well known (cf., for example, \cite{Chernoff73}, \cite[Th.~1.17]{GromLaw83},
\cite{BrMiSh02}) that the operator $\lap_T$ is essentially self-adjoint with initial domain
smooth compactly supported forms. By a slight abuse of notation we will denote by $\lap_T$ also
the self-adjoint extension of this operator to the space $L^2\we$ of square-integrable forms.

\prop{discrete}
For every $T>0$ the spectrum of the operator $\lap_T$ is discrete.
\eprop
We prove the proposition in \refss{prdiscrete} after we study the behaviour of the restriction
of $\lap_T$ to the cylindrical end $\T$.

\subsection{The restriction of the deformed Laplacian to the cylindrical end}\Label{SS:restlap}
Let $v$ denote the vector field on $\tilM$ which corresponds to the 1-form $\ome$ via  the
Riemannian metric $g^{\tilM}$. Since $\ome$ is homogeneous of degree 2, it follows immediately
from \refe{gThomog} that the restriction of $v$ to $\T$ is homogeneous of degree 0, i.e., is
independent of $t$:
\eq{vhomog}
    v(x,st) \ = \tau_{s*}\, v(x,t), \qquad \text{for all} \quad s>1.
\end{equation}

We denote by $\calL_v$ the Lie derivative along $v$. An easy calculations (cf.
\cite[Proposition~11.13]{CFKS}) show that
\begin{equation}\Label{E:DelT=}
    {\lap_T} \ = \
        \frac1{T}\Del+(\calL_v+\calL_v^*)+T|\ome|^2.
\end{equation}
Here $\calL_v^*$ is the formal adjoint of $\calL_v$ and $\Del=\n\n^*+\n^*\n$  stands for the
usual Laplacian on $\tilF$. The operators $(\calL_v+\calL_v^*)$ and $T|\ome|^2$ are symmetric
zero order differential operators, i.e., self-adjoint endomorphisms of the vector bundle
$\Lam^\b(T^*\tilM)\otimes\tilF$ of $\tilF$-valued forms on $\tilM$.

We consider now the restriction of \refe{DelT=} to the cylindrical end $\T$.  For every $x\in
\Gam$, the flat structure defines a trivialization of $\tilF$ (and, hence, of
$\Lam^\b(T^*\tilM)\otimes\tilF$) along the ray $\big\{ (x,t)\in\T:\, t>1\big\}$. Thus for every
$x\in\Gam$ we can consider
\eq{ab}
    a(x,t):= (\calL_v+\calL_v^*), \quad \text{and}\quad b(x,t):= |\ome(x,t)|^2
\end{equation}
as a function of $t$ with values in the space of symmetric matrices.

The proof of \refp{discrete} is based on the following simple lemma, which follows immediately
from \refe{gThomog}, \refe{vhomog}, and the fact that $\ome$ is homogeneous of degree 2 and
non-degenerate at infinity:
\lem{ab-homg}
For all $s>1, \ x\in \Gam, \ t\in [1,\infty)$ the following equalities hold
\eq{ab-homog}
    a(x,st) \ = \ a(x,t), \qquad b(x,st) \ = \ s^2\, b(x,t)> 0.
\end{equation}
\elem

\subsection{Proof of \refp{discrete}}\Label{SS:prdiscrete}
It is well known, cf., for example, \cite[Lemma~6.3]{Shubin99}, that the \refp{discrete} is
equivalent to the following statement: For every  $\eps>0$ there exists a compact set $K\subset
\tilM$ such that, if $\alp\in \Ome^\b(\tilM,\tilF)$ is a smooth compactly supported form, then
\eq{<eps}
    \int_{\tilM\backslash K}\, |\alp|^2\, d\mu \ < \
      \eps \int_\tilW\, \<\lap_T\alp,\alp\>\, d\mu.
\end{equation}
Here, $d\mu$ is the Riemannian volume element on $\tilM$, and $\<\cdot,\cdot\>$ denotes the
Hermitian scalar product on the fibers of $\Lam^\b(T^*\tilM)\otimes\tilF$ induced by the
Riemannian metric on $\tilM$ and the Hermitian metric on $\tilF$.

From \refe{ab-homog}, we conclude that there exists $t_0>1$ such that
\eq{a+b>eps}\notag
    a(x,t) \ + \ b(x,t) \ > \ 1/\eps, \qquad\text{for all} \quad x\in \Gam, \ t\ge t_0.
\end{equation}
Then, using \refe{DelT=} and \refe{ab}, we obtain
\begin{equation*}
    \int_{\Gam\times[t_0,\infty)}\, |\alp|^2\, d\mu \ < \
      \eps\int_{\Gam\times[t_0,\infty)}\, \<(a+b)\alp,\alp\>\, d\mu \ \le \
        \eps\, \int_{\Gam\times[t_0,\infty)}\, \<\lap_T\alp,\alp\>\, d\mu.
\end{equation*}
\hfill$\square$

\th{kerH=}
For every $T>0$, the dimension of the kernel $\Ker\lap_T$ of $\lap_T$ is equal to the dimension
of the relative cohomology $H^\b(\tilM,\tilU^-;\tilF_T)$.

In particular, $\dim\Ker\lap_T^i= \bet_i(\ome,\F)$ for a generic value of $T$.
\eth
The prove of the theorem occupies the next section of the paper.


\section{Computation of the kernel of the Witten Laplacian} \Label{S:prkerH=}

In this section we prove \reft{kerH=}.

\subsection{Complex $\tilOme_T^\b(\tilM, \tilF)$} \Label{SS:OmeTM}
The kernel of the Witten Laplacian $\Ker\lap_T$ has the following cohomological interpretation.
Define
\eq{OmeTM}\notag
    \tilOme^\b_T(\tilM, \tilF) \ = \
     \big\{\,\zet\in L^2\Ome^\b(\tilM, \tilF):\, \nb\zet \in L^2\Ome^\b(\tilM, \tilF)\,  \big\}
\end{equation}
and consider the complex
\begin{equation}\Label{E:com-OmeTM}
 \begin{CD}
    0 \ \to \ \tilOme^0_T(\tilM,\tilF)@>{\nb}>>\tilOme^1_T(\tilM,\tilF)@>{\nb}>>\cdots
            @>{\nb}>>\tilOme^n_T(\tilM,\tilF) \ \to \ 0.
 \end{CD}
\end{equation}
Let $H^\b(\tilOme^\b_T(\tilM,\tilF),\nb)$ denote the cohomology of this complex.

\prop{HodgeT} $H^k(\tilOme^\b_T(\tilM,\tilF),\nb)=\dim \Ker\lap_T^k$\, for all $k=0,1\nek n$.
\eprop
\prf
Since the spectrum of $\lap_T$ is discrete (cf. \refp{discrete}) the space $L^2\Ome^\b(\tilM,
\tilF)$ splits into the orthogonal direct sum of closed subspaces
\eq{Ker+Im}\notag
    L^2\Ome^\b(\tilM, \tilF)= \Ker\lap_T \oplus \mbox{Im}\,\lap_T.
\end{equation}
Hence, every $\varkappa \in \tilOme^\b_T(\tilM, \tilF)$ can be represented as the sum
\eq{DirSum}
     \varkappa  \ = \ \varphi + \lap_T \psi
     \ = \ \varphi + \frac1T\nb^{\phantom{*}}\nb^* \psi + \frac1T\nb^*\nb^{\phantom{*}} \psi,
\end{equation}
where $\varphi\in \Ker\lap_T$ and $\psi\in \Dom(\lap_T)\subset L^2\Ome^\b(\tilM, \tilF)$. Since
$\nb^2=0$ all 3 summands in the right hand side of \refe{DirSum} are mutually orthogonal and
belong to $L^2\Ome^\b(\tilM, \tilF)$.

The ellipticity of $\lap_T$ implies that the forms $\phi$ and $\psi$ are smooth. Hence, so are
the forms $\nb\psi$, $\nb^*\psi$, $\nb^{\phantom{*}}\nb^* \psi$, and
$\nb^*\nb^{\phantom{*}}\psi$. Moreover, $\nb\psi,\, \nb^*\psi\in L^2\Ome^\b(\tilM, \tilF)$
because of the inequality
\eq{SumSq}
    \|\nb \psi \|^2+\| \nb^* \psi\|^2 \ = \ T\<\lap_T\psi, \psi \> \ < \ \infty
\end{equation}
(here $\<\cdot,\cdot \>$ and $\|\cdot\|$ stand for the scalar product and the norm in
$L^2\Ome^\b(\tilM, \tilF)$ respectively). We conclude that $\nb\psi,\, \nb^*\psi\in
\tilOme^\b_T(\tilM, \tilF)$. Hence, \refe{DirSum} implies that the following ``Hodge-type"
decomposition holds
\eq{HodgeT}
    \tilOme^\b_T(\tilM, \tilF) \ = \
        \Ker\lap_T\oplus \IM\,\big(\, \n_T|_{\tilOme^\b_T(\tilM, \tilF)}\,\big)
    \oplus \IM\,\big(\, \n_T^*|_{\tilOme^\b_T(\tilM, \tilF)}\,\big).
\end{equation}
The statement of the proposition follows now from the standard ``Hodge theory" arguments.
\eprf

\reft{kerH=} follows now from the following
\prop{tilcomp}
The cohomology of the complex \refe{com-OmeTM} is isomorphic to the relative cohomology
$H^\b(\tilM,\tilU^-;\tilF_T)$.
\eprop
The proof of the proposition occupies the rest of this section.

\subsection{The cone complex} \Label{SS:ConeCpx}
Let $c$ be an arbitrary positive constant and let $\uc$ be as in \refe{Uc}. We denote by
$j_c:\uc\to\tilM$ the inclusion and by  $j_c^*:\Ome^\b(\tilM,\tilF)\to
\Ome^\b(\uc,\tilF|_{\uc})$ the restriction map.

Consider the cone complex $\Cone_T(j_c^*)$  of the map $j_c^*$ which is defined as follows:
\eq{ConeChains}\notag
    \Cone_T^k (j_c^*) \ =\ \Ome^k (\tilM,\tilF)\oplus \Ome^{k-1}(\uc, \tilF)
\end{equation}
and the differential  $D_{T,c}:\Cone_T^k(j_c^*)\to \Cone_T^{k+1}(j_c^*)$ is given by
\eq{ConeDiff}\notag
    D_{T,c}\,(\eta,\eta_1)\ =\ \Big(\,\nb\eta,\, -\nb\eta_1+ j_c^*\eta\,\Big),\qquad\qquad
    (\eta,\eta_1)\in \Cone_T^\b (j_c^*)
\end{equation}
It is well-known (see eg. \cite{Dold}) that
\[
    H^\b(\Cone_T(j_c^*))\ \cong \ H^\b(\tilM, \uc;\tilF_T).
\]
Thus, to prove \refp{tilcomp} (and, hence, \reft{kerH=}) it suffices to show that the space
$H^\b(\Cone_T(j_c^*))$ is isomorphic to the cohomology $H^\b(\tilOme^\b_T(\tilM,\tilF),\nb)$ of
the complex \refe{com-OmeTM}. This is done in the next subsection, where we construct an
explicit map from $\tilOme_T^\b(\tilM, \tilF)$ to $\Cone_T(j_c^*)$ and show that this map is a
quasi-isomorphism (i.e., induces an isomorphism on cohomology). A similar construction was used
by Farber and Shustin \cite[\S3]{FarberShustin00} in a slightly different situation.

\subsection{A map from $\tilOme^\b_T(\tilM, \tilF)$ to $\Cone_T(j_c^*)$} \Label{SS:OmeT>Cone}
Recall that we assumed that $\ome$ is homogeneous of degree 2 and non-degenerate at infinity.
Thus, cf. \refe{w=df}, there exist $t_0>0$ and a smooth function $f:\Gam\to \RR$ such that
\[
    \ome(x,t) \ = \ d\, \Big(\, f(x)t^2/2\, \Big), \qquad\qquad x\in\Gam, \ t\ge t_0,
\]
and 0 is a regular value of $f$. Changing the parameter $t$ if needed, we can and we will assume
that $t_0=1$. Set
\[
    h(x,t) \ = \ f(x)t^2/2.
\]

The restriction of any form $\zet \in \Ome^\b(\tilM, \tilF)$ to the cylindrical end $\T$ can be
written as
\eq{parperp}
    \zet(x,t) \ = \ dt\wedge \zet_\|(x,t)+\zet_\bot(x,t),\qquad\qquad x\in\Gam, \ \ t>1,
\end{equation}
where $\zet_\|(x,t)$ and $\zet_\bot(x,t)$ are $\tilF$-valued forms on $\Gam$ depending smoothly
on the parameter $t\in (0,\infty)$. Set
\begin{multline}\Label{E:zetbr}
    \breve{\zet}(x,t) \ = \ -e^{-Th(x,t)}\,\int_t^\infty e^{Th(x,\tau)}\zet_\|(x,\tau)\, d\tau
    \\ = \ -\int_t^\infty e^{Tf(x)(\tau^2-t^2)/2}\zet_\|(x,\tau)\, d\tau,
 \qquad  (x,t)\in\uc,
\end{multline}
where we use the flat connection $\n$ to identify the fibers of $\tilF$ along the rays
$\{(x,\tau):\,\tau>0\}$. The integral in \refe{zetbr} converges since $\zet \in
L^2\Ome^\b(\uc,\tilF)$ and $f(x)<0$ for $(x,t)\in \uc$.

Define the map $\Phi_{T,c}:\tilOme_T^\b(\tilM, \tilF) \to \Cone_T(j_c^*)$ by the formula
\eq{Phidef}
    \Phi_{T,c}:\, \zet \ \mapsto \ (\zet, \breve{\zet})
\end{equation}

\lem{ChnMp}
For every $c>0$ the map $\Phi_{T,c}$ is a chain map from $\tilOme_T^\b(\tilM, \tilF; \nb)$ to
$\Cone_T(j_c^*)$, i.e.,
\eq{ChnMp}
    D_{T,c}\circ \Phi_{T,c} \ = \ \Phi_{T,c}\circ \nb.
\end{equation}
\elem
\prf
Let $\zet\in\tilOme_T^\b(\tilM, \tilF)$. The equality \refe{ChnMp} reduces to
\eq{ChnMp2}\notag
    \nb\,\Big[\, e^{-Th(x,t)}\,\int_t^\infty e^{Th(x,\tau)}\zet_\|(x,\tau)\,d\tau\, \Big]+
    \zet(x,t)
    \ = \ - e^{-Th(x,t)}\,\int_t^\infty e^{Th(x,\tau)}(\nb\zet(x,\tau))_\|\,d\tau
\end{equation}
for $(x,t)\in \ucp$. The last equality follows from the following calculation, in which we
denote by $\n_\Gam$ the connection induced by $\n$ on $\Gam$ and by $\n_{\p/\p\tau}$ the
covariant differential along the vector field $\p/\p\tau$:
\begin{multline}\notag
    \nb\,\Big[\, e^{-Th}\,\int_t^\infty e^{Th}\zet_\|\,d\tau\, \Big]+ \zet \ = \
    e^{-Th}\,\n\,\int_t^\infty e^{Th}\zet_\|\,d\tau + \zet
    \\ = \
    -dt \wedge \zet_\| + e^{-Th}\,\int_t^\infty \n_\Gam (e^{Th} \zet_\|) \,d\tau+ \zet
    \ = \
    \zet_\bot + e^{-Th}\,\int_t^\infty \n_\Gam (e^{Th} \zet_\|) \,d\tau
    \\ = \
    e^{-Th}\,\int_t^\infty
       \Big(\, \n_\Gam (e^{Th} \zet_\|) - \n_{\p/\p\tau} (e^{Th} \zet_\bot)\, \Big) \,d\tau
    \ = \
    -e^{-Th}\,\int_t^\infty e^{Th}(\nb\zet)_\|\,d\tau.
\end{multline}
\eprf
Let
\[
    \Phi_{T,c,*}:\, H^\b(\tilOme^\b_T(\tilM,\tilF),\nb)  \ \longrightarrow \
       H^\b(\Cone_T(j_c^*)) \ \cong \ H^\b(\tilM, \uc;\tilF_T)
\]
be the map induced by $\Phi_{T,c}$. To prove \refp{tilcomp} we need now to show that the map
$\Phi_{T,c,*}$ is an isomorphism.

\subsection{Injectivity of $\Phi_{T,c,*}$}
The proof of the injectivity of $\Phi_{T,c,*}$ is based on the following technical lemma.

\lem{alpimproved}
Let $\zet\in \tilOme^\b_T(\tilM,\tilF)$ be a ``harmonic" form, $\lap_T\zet= 0$. Suppose that
\eq{Dalp=zet}
    D_{T,c}\,(\alp,\alp_1)
    \ = \ \Phi_{T,c}\,\zet, \qquad\text{for some}\quad(\alp,\alp_1)\in \Cone_T(j^*_c).
\end{equation}
Then there exist $\eps>0$ and a form $\bet \in \Ome^\b(\tilM,\tilF)$ such that
\begin{enumerate}
\item $\nb\bet = \zet$;
\item  $\bet(x,t) \ \ = \  \breve{\zet}(x,t)\ = \
      -\int_t^\infty e^{\frac12Tf(x)(\tau^2-t^2)}\zet_\| (x,\tau) \, d\tau$,
        \ \ for \  $x\in \cS_\eps$, $t>3$;
\item $\bet_\|(x,t)=0$ \ for \  $x\in \Gam$, $t>3$.
\end{enumerate}
Moreover, any such form $\bet$ belongs to $\tilOme^\b_T(\tilM,\tilF)$.
\elem

\prf
We will construct the form $\bet$ by ``improving" the form $\alp$ in 3 steps.
\subsubsection*{Step 1}
Fix $c^\prime>c$ large enough, so that $\ucp\subset \Gam\times(3,+\infty)$. Choose a smooth
function $\chi_1:\tilM\to [0,1]$ such that $\chi_1\equiv1$ on $\ucp$ and $\chi_1\equiv0$ on
$\tilM\setminus\uc$. Our first approximation of $\bet$ is the form
\eq{alp'}
    \alp^\prime \ = \ \alp-\nb(\chi_1\alp_1).
\end{equation}
From \refe{Dalp=zet} and \refe{alp'} we obtain
\begin{gather}
    \nb \alp^\prime \ = \ \zet \qquad \mbox{on }\ \tilM;\\
     \alp^\prime \ = \ \breve{\zet}   \qquad \mbox{on } \ \ucp.\Label{E:alpp}
\end{gather}
On the cylindrical end $\T$ we can write $\alp^\prime= dt\wedge \alp^\prime_\| +  \alp_\bot'$,
cf. \refe{parperp}. Then the equation \refe{alpp} implies that
\eq{alp'par=0}
    \alp^\prime_\|=0 \qquad\mbox{on }\ \ucp.
\end{equation}

\subsubsection*{Step 2}
Take a new smooth cut-off function $\chi_2:\RR\to[0,1]$, such that $\chi_2(t)= 0$ if $t\le1$ and
$\chi_2(t)= 1$ if $t\ge 2$, and set
\eq{appdef}
    \alp^{\prime\prime}(x,t) \ = \ \alp^\prime(x,t) \ -  \
       \n_T\, \left(\, \chi_2(t)\,e^{-Th(x,t)}\,
        \int_2^t\, e^{Th(x,\tau)}\,\alp^{\prime}_\|(x,\tau)\, d\tau\, \right).
\end{equation}
Again we have $\nb\alp''=\zet$. Since $\chi_2(t)\equiv0$ for $t\ge2$, we also have
\eq{alp''par=0}
    \alp''_\|\ =\ 0\qquad\text{on}\quad\Gam\times(2,\infty),
\end{equation}
so that
\eq{nper=zet}
    \nb\alp''_\bot \ = \ \zet \qquad \text{on}\quad \Gam\times(2,\infty).
\end{equation}
In particular,
\eq{nper=zet2}
    e^{-Th(x,t)}\n_{\p/\p t}\, e^{Th(x,t)}\alp''_\bot(x,t) \ = \ \zet_\|(x,t),
    \qquad\qquad x\in\Gam,\ t>2.
\end{equation}

\subsubsection*{Step 3}
For each $r\in \RR$ denote
\eq{Sr}
    \cS_r \ = \ \big\{\, x\in \Gam:\, f(x)<r/T\, \big\}.
\end{equation}

Recall that we assumed that $\lap_T\zet=0$. A verbatim repetition of the arguments in
\cite[Appendix]{Shubin96Morse} shows that on the cylindrical part $\T$ the form $\zet$ satisfies
the estimate
\eq{estzet}
    |\,\zet(x,t)\,| \ < \ Ce^{-at^2}
\end{equation}
with some positive constants $C$ and $a$. (Here $|\zet(x,t)|$ denotes the norm of $\zet(x,t)$
with respect to the ``conical" Riemannian metric introduced in \refss{metT}). Hence,
\eq{expdecay}\notag
    \big|\, e^{Th(x,t)}\zet_\|(x,t)\big| \ \le \ e^{-(a-r/2)t^2}, \qquad\text{if}\ \ f(x)<r/T,
\end{equation}
i.e., $e^{Th}\zet_\|$ decays exponentially when $t$ tends to infinity and $x\in \cS_{2a}$. Thus
the solution of the ordinary differential equation \refe{nper=zet2} can be written in the form
\eq{alponS-}
    e^{Th(x,t)}\,\alp''_\bot(x,t)
    \ = \ \eta(x) \ - \ \int_t^\infty\, e^{Th(x,\tau)}\,\zet_\|(x,\tau)\,d\tau,
    \qquad\qquad x\in \cS_{2a},\ t>2,
\end{equation}
where $\eta(x)$ is a differential form on $S_{2a}\subset \Gam$ (independent of $t$).

From \refe{nper=zet} and \refe{alponS-} we see that on $\cS_{2a}\times(2,\infty)$
\eq{zetbot}
    \zet_\bot \ = \ e^{-Th}\,\n_\Gam\, e^{Th}\,\alp''_\bot \ = \
    e^{-Th}\,\n_\Gam\,\eta \ - \ e^{-Th}\, \int_t^\infty\, \n_\Gam\, e^{Th}\,\zet_\|\, d\tau.
\end{equation}
Since $\nb\zet=0$ we have
\eq{nTGam}\notag
    0 \ = \ \big(\, \nb\zet\, )_\| \ = \ -\, e^{-Th}\,\n_\Gam\, e^{Th}\,\zet_\|
    \ + \ e^{-Th}\,\n_{\p/\p t}\, e^{Th}\,\zet_\bot,
\end{equation}
i.e.,
\eq{nG=nt}\notag
    \n_\Gam\, e^{Th}\,\zet_\| \ = \ \n_{\p/\p t}\, e^{Th}\,\zet_\bot.
\end{equation}
Substituting this equality into \refe{zetbot}, we see that on $\cS_{2a}\times(2,\infty)$
\eq{zetbot2}\notag
    \zet_\bot \ = \
    e^{-Th}\,\n_\Gam\,\eta \ - \ e^{-Th}\, \int_t^\infty\, \n_{\p/\p\tau}\, e^{Th}\,\zet_\bot\, d\tau \ = \
    e^{-Th}\,\n_\Gam\,\eta \ + \ \zet_\bot
\end{equation}
Hence,
\eq{nbet}\notag
    \n_\Gam\, \eta(x) \ = \ 0, \qquad \text{if} \quad x\in \cS_{2a},
\end{equation}
i.e., $\eta$ is a closed form on $\cS_{2a}$.

Set $\gam= \int_2^\infty{}e^{Th}\alp'_{\|}d\tau$. Comparing \refe{alponS-} with \refe{appdef}
and using \refe{alpp} we see that $\eta= \n\gam$ on $\ucp$. Since, by the construction, $\gam$
is independent of $t$, we obtain $\eta=\n_\Gam\gam$ on $\ucp$. But, for all $x\in \cS_0$, there
exists $t$ such that $(x,t)\in \ucp$. Therefore, $\eta= \n_\Gam\gam$ on $\cS_0$. In particular,
the restriction of $\eta$ to $\cS_0$ represents the zero class in $H^\b(\cS_0)$.

Since $0$ is a regular value of $f$ there exists $\eps\in (0,a)$ such that $\cS_0$ is a
deformation retract of $\cS_{2\eps}$. Hence $\eta$ also represents the trivial class in
$H^\b(S_{2\eps})$, i.e., there exists a form $\xi\in \Ome^\b(\cS_{2\eps})$ such that
\eq{xi}
    \n_\Gam\,\xi(x) \ = \ \eta(x), \qquad x\in \cS_{2\eps}.
\end{equation}
Consider a smooth cut-off function $\chi_3:\Gam\to[0,1]$, such that $\chi_3 = 1$ on
$\cS_{\frac43\eps}$ and $\chi_3 = 0$ on $\Gam\backslash \cS_{\frac53\eps}$. Finally, set
\eq{bet2}
    \bet \ = \ \alp'' -  \nb\, \Big(\, e^{-Th(x,t)}\,\chi_3(x)\,\chi_2(t)\,\xi\,\Big).
\end{equation}
Since $\bet$ is obtained by subtracting exact forms from $\alp$, condition (1) of the lemma is
satisfied. From \refe{alp''par=0}, \refe{xi}, \refe{bet2}, and \refe{alponS-} we obtain
\eq{alp'''S}\notag
    \bet(x,t) \ = \ \alp''_\bot(x,t) \ - \ e^{-Th(x,t)}\,\eta(x) \ = \
     - \int_t^\infty\, e^{\frac12Tf(x)(\tau^2-t^2)}\,\zet_\|(x,\tau)\,d\tau,
     \qquad\text{for}\quad x\in \cS_\eps, t>3.
\end{equation}
Thus condition (2) of the lemma is satisfied.

On $\Gam\times(3,+\infty)$ form $\bet$ is obtained from $\alp''$, whose parallel component
vanishes, by subtracting $\nb\, (\, e^{-Th(x,t)}\,\chi_3(x)\,\,\xi\,)$, which also does not have
a parallel component. Therefore, condition (3) of the lemma is also satisfied.

It remains to prove now that $\bet$ is in $\tilOme^\b_T(\tilM, \tilF)$. The estimate
\refe{estzet} implies that $\nb\bet=\zet$ is in $L^2\Ome^\b(\tilM,\tilF)$, so we only need to
show that $\bet\in L^2\Ome(\tilM,\tilF)$. Since the form $\bet$ is smooth, its restriction to
every compact set is square integrable. Hence, it is enough to show that the restrictions of
$\bet$ to $\cS_\eps\times(3,\infty)$ and $(\Gam\backslash{}\cS_\eps)\times(3,\infty)$ are square
integrable.

\subsubsection*{Restriction of $\bet$ to $\cS_\eps\times(3,\infty)$}
Using the estimate \refe{estzet}, we see that on $\cS_\eps\times(3,\infty)$
\eq{onAeps}\notag
    |\,\bet(x,t)\,|
    \ = \ \left|\,\int_t^\infty\,
        e^{\frac12Tf(x)(\tau^2-t^2)}\,\zet_\| (x,\tau)\, d\tau \,\right|
    \ \le \ C\int_t^\infty\, e^{-(a-\eps/2)\tau^2}\, d\tau.
\end{equation}
Since $\eps<a$ we conclude that $\bet$ is decaying exponentially with $t$ and, hence, is square
integrable on $\cS_\eps\times(3,\infty)$.

\subsubsection*{Restriction of $\bet$ to $(\Gam\backslash S_\eps)\times(3,\infty)$}
Since $\bet_\|=0$ on $\Gam\times (3,\infty)$, taking the parallel component of the equality
$\zet= \n_T\bet$, we obtain
\eq{restbet}\notag
 e^{-Th(x,t)}\,\n_{\p/\p t}\, e^{Th(x,t)}\, \bet(x,t) \ = \ \zet_\|(x,t).
\end{equation}
Solving this ordinary differential equation we get
\eq{alp=int}\notag
    e^{Th(x,t)}\,\bet(x,t) \ = \ \bet(x,3)
    \ + \ \int_3^t\, e^{Th(x,\tau)}\, \zet_\|(x,\tau)\, d\tau.
\end{equation}
Therefore, setting  $B=\max\limits_{x\in \Gam} |\bet(x,3) |$ and using \refe{estzet}, we see
that on $(\Gam\backslash S_\eps)\times(3,\infty)$ the following estimate holds
\begin{multline}\notag
    \left|\,\bet(x,t)\, \right|
    \ \le \ B\,e^{-\frac12Tf(x)t^2}
        + \int_3^t\, e^{\frac12 Tf(x)(\tau^2-t^2)}\,|\zet_\|(x,\tau)|\, d\tau
    \\ \le \
    B\, e^{-\eps t^2/2} \ + \
    \int_3^t\, e^{\frac12 \eps(\tau^2-t^2)}\,e^{-a\tau^2}\, d\tau
    \ \le \
    B\, e^{-\eps t^2/2} \ + \
     e^{-\frac12\eps t^2}\,\int_3^t\, e^{-(a-\frac12 \eps)\tau^2}\, d\tau.
\end{multline}
Thus $\bet(x,t)$ decays exponentially when $t$ tends to infinity and, hence, is square
integrable on $(\Gam\backslash S_\eps)\times(3,\infty)$. The proof of the lemma is now complete.
\eprf

\cor{Phimono} The map $\Phi_{T,c,*}$ is injective.
\ecor

\prf
Let $\Upsilon\in H^\b(\tilOme^\b_T(\tilM,\tilF))$ be a cohomology class, whose image under
$\Phi_{T,c,*}$ is trivial. By \refp{HodgeT}, the representative $\zet\in \Upsilon$ can be chosen
in such a way that $\lap_T\zet=0$. Then  \refl{alpimproved} implies that there exist
$\bet\in\tilOme^\b_T(\tilM,\tilF)$, such that  $\nb\bet= \zet$. Hence,  $\Upsilon= 0$.
\eprf

\subsection{Surjectivity of $\Phi_{T,c,*}$}
In order to prove that $\Phi_{T,c,*}$ is surjective we need the following auxiliary lemma,
which, roughly speaking, states that every relative cohomology class has a bounded
representative.

Recall that the sets $\cS_r\subset \Gam$, were defined in \refe{Sr}. Since $0$ is a regular
value of $f$ we can take $\eps>0$ small enough, such that for all $\eps'\in[0,\eps]$, the set
$\cS_{\eps'}$ is homotopically equivalent to $\cS_0= U^-$. For any $\del>1$ denote by
$j_{\del,\eps}$ the inclusion map of $\cS_\eps\times[\del,+\infty)$ into $\tilM$. Let
$\Cone_T(j_{\del,\eps}^*)$ be the complex
\[
    \Cone_T^k(j_{\del,\eps}^*)\ =\
    \Ome^{k}(\tilM,\tilF)\oplus\Ome^{k-1}(\cS_\eps\times[\del,+\infty)),\qquad
    k=0\nek\dim \tilM
\]
with the differential
\eq{DTde}\notag
    D_{T,\del,\eps}(\eta,\eta_1)\ =\
       \big(\,\nb\eta, -\nb\eta_1+ j_c^*\eta\,\big),\qquad\qquad
    (\eta,\eta_1)\in \Cone_T^\b (j_c^*).
\end{equation}
Also for any $\del>0$ denote by $\tilM_\del$ the compact set $M\cup(\Gam\times[1,\del))\subset
\tilM$.

\lem{aapbdd}
Every cohomological class in $H^\b(\Cone_T(j^*_{\del,\eps}))$ has a representative
$(\alp,\alp_1)$, such that the restrictions of the forms $e^{Th}\alp$, $e^{Th}\alp_1$, $\nb
e^{Th}\alp$, and $\nb e^{Th}\alp_1$ to the cylindrical end $\calT$ are bounded with respect to
the Riemannian metric introduced in \refss{metT}.
\elem

\prf
As in the proof of \refl{alpimproved} we start with an  arbitrary representative
$(\alp',\alp'_1)$ and ``improve it''.

For all $\eps'\in (0,\eps)$ and $\del'>\del$, the pairs $(\tilM,\cS_\eps\times[\del,+\infty))$
and $(M_{\del'}, \cS_{\eps'}\times[\del,\del'))$ are homotopically equivalent. Moreover, there
exists a smooth map $\calR:\tilM\to\tilM_{\del'}$, such that:
\begin{itemize}
\item[(i)] $\calR$ is a homotopy equivalence of the topological pairs $(\tilM,\cS_\eps\times[\del,+\infty))$ and $(\tilM_{\del'}, \cS_{\eps'}\times[\del,\del'))$,
\item[(ii)] the restriction of $\calR$ to $\tilM_{\del}$ is the identity operator,
\item[(iii)] the tangent map $T\calR$ is bounded (with respect to the Riemannian
metric introduced in \refss{metT}).
\end{itemize}
On  the cylindrical end $\calT=\Gam\times[1,+\infty)$ the differentials $\nb$ and
$D_{T,\del,\eps}$ can be written as follows
\begin{equation}\notag
    \nb \ = \ e^{-Th}\,\nb\, e^{Th}, \qquad
    D_{T,\del,\eps}\ =\  e^{-Th}\,D_{0,\del,\eps}\,e^{Th}.
\end{equation}
Properties (i) and (ii) of the map $\calR$ imply that on $\calT$
\eq{mumu1}
    e^{Th} (\alp',\alp'_1)\ -\ \calR^*e^{Th}(\alp',\alp'_1)
    \ = \ D_{0,\del,\eps} (\mu,\mu_1)
\end{equation}
for some $\mu\in \Ome^\b(\calT, \tilF)$ and  $\mu_1 \in\Ome^\b(\calS_\eps\times[\del,+\infty))$.
Let $\chi:\RR\to [0,1]$ be a smooth nondecreasing function, such that
\[
    \chi(t)=0 \quad \text{for } \
    t\le\scriptstyle\frac{1+\del}{4},\displaystyle,\quad \text{and}\quad
    \chi(t)=1 \quad \text{for } \
    t\ge\scriptstyle\frac{1+\del}{2}\displaystyle.
\]
We also denote by $\chi$ its natural extension to $\tilM$. Set
\eq{alpalp1}
    (\alp,\alp_1) \ = \
    (\alp',\alp'_1)\ - \
    D_{T,\del,\eps}\,\left(\,e^{-Th}\chi\cdot(\mu,\mu_1)\,\right).
\end{equation}
Then $(\alp,\alp_1)$ and $(\alp',\alp'_1)$ represent the same cohomological class in
$(\Cone_T(j^*_{\del,\eps}))$. By construction, $\chi(t)=1$ for $t\ge \frac{1+\del}{2}$. Hence,
using \refe{mumu1} and \refe{alpalp1}, we obtain
\begin{multline}
    e^{Th}\, (\alp,\alp_1)\ =\
    e^{Th}\,\Big(\,(\alp',\alp_1') \ - \
        D_{T,\del,\eps}\,\left(e^{-Th}(\mu,\mu_1)\right)\,\Big)\\
    \ = \ e^{Th}\,(\alp',\alp_1')\ -\ D_{0,\del,\eps}(\mu,\mu_1)
    \ = \
 \calR^*\,e^{Th}\,(\alp',\alp'_1),
  \qquad\text{on }\quad\Gam\times[2,+\infty).\Label{E:R2}
\end{multline}
The pair $\calR^*e^{Th}(\alp',\alp_1')$ is the pull--back of the restriction of
$e^{Th}(\alp',\alp_1')$ to the finite part $\tilM_{\del'}$ of $\tilM$. Since $T\calR$ is
bounded, \refe{R2} implies that $e^{Th}\alp$ and $e^{Th}\alp_1$ are also bounded together with
their differentials.
\footnote{Recall that we use the conical metric on $\calT$ introduced in \refss{metT}. Thus, if,
for example, we consider a form $\gam\in \Ome^\b(\Gam)$ as a constant form on $\calT=
\Gam\times(0,\infty)$ then the norm $|\gam|_{(x,t)}$ of $\gam$ at the point $(x,t)\in \calT$
decays as $1/t$ as $t\to \infty$.}
\eprf

\cor{Phiepi} The map $\Phi_{T,c,*}$ is surjective.
\ecor

\prf
Since the topological pairs $(\tilM,\uc)$ are homotopically equivalent for different $c>0$, it
suffices to prove the corollary for sufficiently large $c$. Therefore, we can assume that $c$ is
large enough, so that $\uc\subset\Gam\times[3,+\infty)$.

Let $\Upsilon\in H^\b(\Cone_T(j_c^*))$ be an arbitrary cohomology class. Take $\eps>0$ small
enough, so that $U^-=\cS_0$ is a deformation retract of $\cS_{\eps'}$ for all
$\eps'\in[0,2\eps]$.  Let
\eq{R:}\notag
    R:\, \Cone_T(j^*_{2,2\eps})\ni (\alp,\alp_1)\ \mapsto \
    (\alp,\alp_1|_{\uc})\in \Cone_T (j_c^*)
\end{equation}
be the restriction map. Since $\uc$ is a deformation retract of $\cS_{2\eps}\times[2,+\infty)$,
the map $R$ induces an isomorphism of cohomology. Therefore, $\Upsilon$ can be represented by a
pair $R(\alp,\alp_1)$, where $(\alp,\alp_1)\in \Cone_T(j_{2,2\eps}^*)$. By \refl{aapbdd}, the
pair $(\alp,\alp_1)$ can be chosen in such a way that the forms  $e^{Th}\alp$,
$e^{Th}\alp_1$,$e^{Th}\nb\alp$, and $e^{Th}\nb\alp_1$ are all bounded.

Let $\chi: \tilM\to[0,1]$ be a smooth function, such that
\begin{equation}\notag
    \chi=0 \quad \text{on } \ \tilM\setminus (\cS_{2\eps}\times [2,+\infty)), \qquad
    \text{and}\qquad
    \chi=1  \quad \text{on } \ \cS_{\eps}\times [3,+\infty).
\end{equation}
and $d\chi$ is bounded on $\calT$ with respect to the Riemannian metric introduced in
\refss{metT}. Set
\eq{bet=alp-}\notag
    \bet \ =\ \alp\ -\ \nb (\chi\alp_1).
\end{equation}
Then the pair $(\bet,0)$ is cohomologous to the pair $(\alp,\alp_1)$ in the $\Cone_T(j_c^*)$. In
addition, $\bet$ vanishes on $\cS_\eps\times[3,+\infty)\supset \uc$, hence
$[\Phi_{T,c}\bet]=\Upsilon$. In order to finish the proof of the corollary it remains to show
that $\bet\in\tilOme^\b_T(\tilM,\tilF)$. Since $\nb\bet=\nb\alp=0$, we only need to check that
$\bet$ is square integrable on the cylindrical end $\calT$. Since $\bet$ vanishes on
$\calS_\eps\times[3,+\infty)$ we only need to consider its restriction to
$(\Gam\setminus\calS_\eps)\times[3,+\infty)$. The boudedness of $e^{Th}\alp$, $e^{Th}\alp_1$ and
$d\chi$ implies that $|e^{Th}\bet|<B$ for some $B>0$. Hence, on
$(\Gam\setminus\calS_\eps)\times[3,+\infty)$ we have
\eq{|bet|}\notag
    |\bet(x,t)| < e^{-Th(x,t)}B=Be^{-\frac12f(x)t^2}<Be^{-\frac12\eps{}t^2},
\end{equation}
hence $\bet$ is square integrable and the map $\Phi_{T,c,*}$ is surjective.
\eprf

The proof of \reft{kerH=} is complete. \hfill$\square$


\section{Deformation of the Laplacian on the normal bundle to a minimizing manifold}
\Label{S:witN}

In this section for each critical subset $C$ of $\ome$ we introduce a manifold $\tilN_C$ with a
cylindrical end, whose compact part may be identified with a neighborhood of $C$ in $\tilM$. We
also introduce a family of connections $\nb^\nce$ on the lift of $\tilF$ to $\nce$, which is
similar to the Novikov deformation \refe{nab-t}. We apply the machinery of \refs{witten} to the
Laplacians $\lap_T^{\N_C}$ of the connections $\nb^\nce$. In particular, we show that the
spectrum of $\lap_T^{\N_C}$ is discrete and does not depend on $T$ and we compute the kernel of
$\lap_T^{\N_C}$. The dimension of the kernel $\lap_T^{\N_C}$ encodes the information about the
topology of the critical subset $C$, cf. \reft{lapE}.

The significance of the operators $\lap_T^{\N_C}$ is that their direct sum
$\bigoplus_{C\in\bfC}\lap_T^{\N_C}$ plays the role of a {\em model operator} (cf.
\cite{Shubin96Morse}) for the family $\lap_T$ of Laplacian on the manifold $\tilM$. More
precisely, this statement is formulated in the next section (\reft{model}).

\subsection{Construction of the manifold $\tilN_C$} \Label{SS:tilN}
Let $C\subset \bfC$ be a critical subset of \w, cf. \refss{nondeg}. Let $h_C,\, \Sig_C$ be as in
\refd{kirwan} and let $\nu(\Sig_C), \, \calW\subset \nu(\Sig_C)$ be as in \refr{C3}.

Fix a Hermitian metric $h^{\nu(\Sig_C)}$ and a Hermitian connection $\n^{\nu(\Sig_C)}$ on the
vector bundle $p:\nu(\Sig_C)\to\Sig_C$.

Fix a small enough number $\eps>0$ so that
\begin{itemize}
\item
 \(\displaystyle
    \big\{y\in \nu(\Sig_C):\, h_C(p(y))+\frac{|y|^2}2 \le (7\eps)^2\big\}
    \ \subset \ \calW.
  \)
\item
The restriction of $h_C$ to $\Sig_C\cap\calW$ does not have critical values on the interval
$(1,7\eps)$.
\end{itemize}

Set
\[
    V_{C,r} \ := \ \big\{\, x\in \Sig_C:\, h_C(x)<r^2\, \big\}.
\]
Using the gradient flow of $h_C|_{\Sig_C}$ with respect to some arbitrary metric, we can
construct a family of smooth maps $m_t:V_{C,7\eps}\to V_{C,7\eps}$ $(0<t\le 1)$ such that
$m_{t_1}\circ{}m_{t_2}= m_{t_1t_2}$ and
\eq{mt}
  \begin{aligned}
    h_C(m_t(x)) \ &= \
        t^2\, h_C(x), \qquad x\in V_{C,7\eps}\backslash V_{C,\eps/2}, \ 1/2<t\le 1,\\
    h_C(x) \ \ge \ h_C(m_t(x)) \ &\ge \
        t^2\, h_C(x), \qquad x\in  V_{C,7\eps}, \ 0<t\le 1.
  \end{aligned}
\end{equation}

For each $y\in p^{-1}V_{C,7\eps}$ let
\[
    \tilm_t(y) \ \in \ p^{-1}\big(\,m_t(p(x))\,\big)
\]
denote the horizontal lift of the curve $t\mapsto m_t(p(y))$. Then $|\tilm_t(y)|=|y|$ for all
$t\in (0,1]$.

Set
\eq{BrSr1}
  \begin{aligned}
    S_{C,1} \ &= \
      \big\{\, y\in p^{-1}V_{C,7\eps}:\, h_C\big(p(y)\big)+|y|^2/2=\eps^2\,\big\},\\
    B_{C,1} \ &= \
    \big\{\, y\in p^{-1}V_{C,7\eps}:\, h_C\big(p(y)\big)+|y|^2/2\le \eps^2\, \big\}.
  \end{aligned}
\end{equation}

We will now construct a neighborhood $\calU\subset \calW$ of $B_{C,1}$ and a diffeomorphism
\[
    \phi:\, S_{C,1}\times(1,7) \ \longrightarrow \ \calU\backslash{}B_{C,1}
\]
such that, if $\phi(y_1,t)= y, \ y_1\in S_{C,1}, \ t\in (1,7)$, then
\begin{align}
    m_{1/t}\big(\,p(y)\,\big) \ &= \ p(y_1); \Label{E:py=py1}\\
    h_C\big(\,i(y)\,\big) \ &= \ t^2\, h_C\big(\,i(y_1)\,\big). \Label{E:y=ty1}
\end{align}
The equation \refe{py=py1} implies that there exists $\tau=\tau(y)$ such that
\eq{y=tauy1}
    y \ = \ \tau\, \tilm_t(y_1).
\end{equation}
Then $|y|^2= \tau^2 |y_1|^2$. The equation \refe{y=ty1} is equivalent now to the equation
\eq{eqtau}\notag
    h_C\big(\, m_t(p(y_1))\,\big) \ - \ \tau^2\,\frac{|y_1|^2}2 \ = \
    t^2\, \Big[\, h_C\big(\, p(y_1)\,\big) \ - \ \frac{|y_1|^2}2\,\Big].
\end{equation}
or
\eq{eqtau2}
    \big(\,t^2-\tau^2\,\big)\, \frac{|y_1|^2}2
    \ \le \
    t^2\, h_C\big(\,p(y_1)\,\big) \ - \ h_C\big(\,m_t(p(y_1))\,\big).
\end{equation}
To construct the diffeomorphism $\phi$ it is enough now to show that the equation \refe{eqtau2}
has a unique positive solution $\tau= \tau(y)>0$.

Let us, first, consider the case $h_C(i(y_1))\ge 0$. Then $h_C(p(y_1))\ge \eps^2/2$. In
particular, $p(y_1)\in V_{C,7\eps}\backslash{}V_{C,\eps/2}$. Then, by \refe{mt}, the right hand
side of \refe{eqtau2} vanishes. Hence, we can set $\tau=t$.

Consider now the case $h_C(i(y_1))< 0$. Then $|y_1|/2 > h_C(p(y_1))$. Hence, \refe{eqtau2} and
\refe{mt} imply that
\[
    0 \ \le \ t^2\, h_C\big(\,p(y_1)\,\big) \ - \ h_C\big(\,m_t(p(y_1))\,\big) \ < \
    t^2\, \frac{|y_1|^2}2.
\]
Thus there exists a unique $\tau\in (0,t]$ satisfying \refe{eqtau2}.

Summarizing the above 2 cases, we see that there always exists a unique $\tau\in (0,t]$ which
satisfy \refe{eqtau2}. We now use this $\tau$ to define $y=\phi(y_1,t)$ by \refe{y=tauy1}.

In the sequel we will identify the pair $(y_1,t)$ and the vector $y=\phi(y_1,t)$ and we will
write $y= (y_1,t)$.

For $r\in [1,7)$ set
\eq{BrSr}
  \begin{aligned}
    S_{C,r} \ &= \
    \big\{\, y\in \calW:\, y=(y_1,r), \ y_1\in S_{C,1} \big\}, \\
    B_{C,r} \ &= \
    B_{C,1}\cup\big\{\, y\in \calW:\, y=(y_1,t), \ y_1\in S_{C,1}, \ 1\le t\le r\,\big\}.
  \end{aligned}
\end{equation}
Then $B_{C,6}$ is naturally diffeomorphic to the union $B_{C,1}\cup\big(\,
S_{C,1}\times[1,6]\,\big)$. Hence, we can define a manifold  with a cylindrical end
\eq{tilN}\notag
    \nce \ := \ B_{C,6}\cup \big(\, S_{C,1}\times [1,\infty)\, \big)
    \ = \ B_{C,6} \sqcup \big(\, S_{C,1}\times (6,\infty)\, \big).
\end{equation}
It follows from \refe{y=ty1} that there exists a smooth function $\tilh_C:\nce\to \RR$ such that
\eq{tilhC}\notag
    \tilh_C(y) \ = \
     \begin{cases}
        h_C\big(\, i(y)\,\big), \quad&\text{for} \quad y\in B_{C,6};\\
        t^2h_C\big(\,i(y_1)\,\big), \quad&\text{for}
        \quad y=(y_1,t)\in S_{C,1}\times[1,\infty).
    \end{cases}
\end{equation}

\subsection{Deformation of the covariant differential} \Label{SS:nbE}
Choose a smooth non-decreasing function $\sig:[0,+\infty) \to [0,1]$, such that
\[
    \sig(t) = 0 \quad\text{for } \ t\in [0,2], \qquad \text{and}\qquad
    \sig(t)= 1 \quad\text{for } \ t\in [3,+\infty).
\]
For each $s>0$, set
\eq{phi(s,y)}\notag
    \phi(s,y) \ =
    \begin{cases}
        1, \qquad & \text{if}\quad y\in B_{C,2}\\
            s^{\sig(t)}, \qquad & \text{if} \quad y=(x,t)\in S_{C,1}\times(2,\infty)
    \end{cases}
\end{equation}
Since $\sig$ vanishes on $[0,2]$, the function $\phi(s,y)$ is well-defined and smooth on $\nce$.
The vector bundle $\tilF$ induces naturally a flat vector bundle on $\nce$, which we will also
denote by $\tilF$. Let $\nabla^\nce$ denote the covariant differential operator determined by
the flat structure on $\Ff$. (We use the superscript $\nce$ here, to distinguish this operator
from the covariant differential on the bundle $\F$ over the manifold $\tilM$). Consider the
one-parameter deformation $\nb^\nce$ of this covariant differential defined by the formula
\eq{dE}\notag
    \nb^\nce\  \ =\
     e^{-\phi({T},y)\cdot \tilh_C(y)}\,\n^\nce\, e^{\phi({T},y)\cdot \tilh_C(y)}
\end{equation}
\rem{phiisOK}
Note that, for each $T>0$, the equation \refe{dE} is equivalent to \refe{nab-t} with
\[
    \ome\ =\ \frac1T\,d\,\big(\,\phi({T},y)\,\tilh_C\,\big).
\]
Therefore, we can apply the results of the previous sections of the paper to the connection
$\nb^\nce$. Moreover, since the form $\ome$ is exact on $\nce$ the Novikov numbers will be
everywhere replaced by the usual Betti numbers.
\erem

For every $s>0$, let $r_s:\nce\to \nce$ be the map defined by
\eq{phist}\notag
    r_s(y) \ = \
    \begin{cases}
     y, \qquad & \text{if}\quad y\in B_{C,2},\\
     \big(\,x,\,\phi(s,y)t\,\big), \qquad & \text{if} \quad y=(x,t)\in S_{C,1}\times(2,\infty).
    \end{cases}
\end{equation}
If $y\in \nce$, we denote by $\Ff_y$ the fiber of $\Ff$ over $y$. The flat connection on $\Ff$
gives a natural identification of the fibers $\Ff_y$ and $\Ff_{sy}$.  Hence, the map
$r_s:\nce\to \nce$ defines the ``pull-back'' map
\[
    r^*_s:\,\Ome^\b(\nce,\Ff)\ \to\ \Ome^\b(\nce,\Ff).
\]
\lem{de-de'} For every $T>0$
  \begin{equation}\Label{E:dE'}
    \nb^\nce \ = \ \r\, e^{-\tilh_C}\,\n^\nce\, e^{\tilh_C}\,(\r)^{-1}.
  \end{equation}
\elem

\begin{proof}
Since $\r$ commutes with $\n^\nce$, we can write
\begin{multline}\notag
    \r\, e^{-\tilh_C}\, \n\, e^{\tilh_C}\, (\r)^{-1}
      \ =\ \r\, e^{-\tilh_C}\,(\r)^{-1}\, \nb\, \r\, e^{\tilh_C}\, (\r)^{-1}\\
    \ = \ e^{-\phi({T},y) \tilh_C}\,\n^\nce\, e^{\phi({T},y) \tilh_C}
    \ = \ \nb^\nce.
\end{multline}
Here the last line follows from the identity
\begin{equation}\Label{E:tft}
    \r\, \tilh_C(y)\, (\r)^{-1} \ =\  \phi({T},y)\, \tilh_C(y),
  \end{equation}
where $\tilh_C$ is identified with the operator of multiplication by $\tilh_C$.
\end{proof}
\subsection{Deformation of the Laplacian} \Label{SS:witN}
Let us choose a Riemannian metric $g^\nce_1$ on $\nce$ and a Hermitian metric $h_1^{\Ff}$ on
$\Ff$ which on the cylindrical end $S_{C,1}\times(1,\infty)$ have the form described in
\refss{metT}. Consider the families of metrics $g^\nce_T$ and $h^{\Ff}_T$ depending on $T$ given
by
\eq{hgt}
    g^\nce_T\ =\ \frac1T\,\r{}\,g_1^\nce, \qquad\qquad h^{\Ff}_T \ = \ \r\, h^{\Ff}_1
\end{equation}

\defe{laplE}
The {\em Witten Laplacian} $\lap^\nce$ of the bundle $\nce$ associated to the metrics $g_T^\nce,
h^{\Ff}$ is defined by the formula
\begin{equation}\Label{E:lapl}
    \lap^\nce_T\ = \ \frac1T\,\big(\,\nb^\nce\,\nb^{\nce*}\,+\,\nb^{\nce*}\,\nb^\nce\,\big),
\end{equation}
where $\nb^{\nce*}$ denote the formal adjoint of $\nb^\nce$ with respect to the metrics
$g^\nce_T,\ h^{\Ff}_T$. We denote by $\lap^{\nce,p}_T$  the restriction of $\lap^{\nce}_T$ to
the space of $p$-forms.
\edefe

\prop{del=del}For every  $T>0$,
  \begin{equation}\Label{E:del=del}
    \lap^\nce_T \ = \ \r\, \lap^\nce_1\, (\r)^{-1}.
  \end{equation}

  In particular, the operators $\lap^\nce_T$ and $\lap^\nce_1$
  have the same spectrum.
\eprop
\begin{proof}
\refl{de-de'} implies that
\eq{nT1}
    \n^\nce_T \ = \ \r \,\n_1^\nce\, (\r)^{-1}.
\end{equation}
The metrics  $g^\nce_T=\frac1T\r{}g_1^\nce$ and $h^{\Ff}_T=\r h^{\Ff}_1$ are defined in such a
way that the map
\[
    r_{\sqrt{T}}:\, (\nce, Tg_T^\nce)\ \to\ (\nce, g_1^\nce)
\]
is an isometry (here $(\nce, Tg_T^\nce)$ and $(\nce, g_1^\nce)$ denote $\nce$ considered as a
Riemannian manifold with the Riemannian metric $Tg_T^\nce$ and $g_1^\nce$ respectively).
Moreover, the pull-back map
\[
    \r:\, L^2\big(\,\Ome^\b(\nce,\tilF), g_1^\nce\,\big)\ \to\
        L^2\big(\,\Ome^\b(\nce,\tilF), Tg_T^\nce\,\big)
\]
is also an isometry (here $L^2(\Ome^\b(\nce,\tilF),g)$ denotes the space of differential forms,
which are square-integrable with respect to the scalar product induced by the Riemannian metric
$g$).

We conclude that the adjoint operator to the differential $\n$ is equal to
\eq{nstar}
     \nb^* \ = \ \frac1T\,\r\, \n_1^*\, (\r)^{-1}.
\end{equation}
The equality \refe{del=del} follows now from \refe{nT1} and \refe{nstar}.
\end{proof}

Recall from \refss{calM} that for each critical subset $C$, we denote by $\ind(C)$ the dimension of
the fibers of the bundle $\nu(\Sig_C)\to \Sig_C$ and by $o(C)$ the orientation bundle of
$\nu(\Sig_C)$, considered as a flat line bundle.

Recall that $\eps>0$ was chosen in \refss{tilN}. The following theorem gives the main  spectral
properties of the Witten Laplacian $\Del^{\tilN_C}_T$.
\th{lapE} \ (1) \ For each $T>0$ the operator $\lap_T^\nce$ has a discreet spectrum. More
    precisely, there exists an orthonormal basis of {\em smooth}
    eigenforms with eigenvalues
    $\lam_j\in{}\RR$, such that $\lam_j\to{}\infty$ as $j\to{}\infty$.
    The spectrum of $\lap_T^\nce$ coincides with the set of
    all eigenvalues $\{\lam_j\}$.

(2) \ Let $\lap^{\nce,p}_T$ denote the restriction of $\lap^\nce_T$ on the
    space of $p$-forms. Then
    \eq{kerlapnce}
       \dim\Ker \lap^{\nce,p}_T \ = \
        \dim {H}^{p-\ind(C)}
           (V_{C,\eps},\F_{|_{V_{C,\eps}}}\otimes o(\Sig_C)_{|_{V_{C,\eps}}})
    \end{equation}
    for any $p=0,1,\dots, \dim M$.
\eth
\begin{proof}
The part (1) of the theorem is a particular case of \refp{discrete}.

Since the deformation $\nb^\nce$ is defined using the exact form $\frac1Td(\phi{}\tilh_C)$ (cf.
\refr{phiisOK}), \reft{kerH=} implies that
\eq{kerlapnce2}
       \dim\Ker \lap^{\nce,p}_T \ = \  \dim H^{p}(\nce,\tilU_{C,c}^-;\F),
\end{equation}
where $c>0$ and
\[
    \tilU_{C,c}^- \ = \ \big\{\, y\in \nce:\, \tilh_C(y)<c\,\big\}.
\]
Set $U_{C}^-:= \{y\in S_{C,1}:\, \tilh_C(y)<0\}$. Clearly, the topological pair
$(\nce,\tilU_{C,c}^-)$ is homotopic to the pair $(B_{C,1},U_{C}^-)$.

Let $\nu_{C,\eps}= p^{-1}V_{C,\eps}$ be the restriction of the bundle $p:\nu(\Sig_C)\to \Sig_C$
to $V_{C,\eps}$. Since the pair $(B_{C,1},U_{C}^-)$ is homotopic to the pair
$(\nu_{C,\eps},\nu_{C,\eps}\backslash{V_{C,\eps}})$, the right hand side of \refe{kerlapnce2} is
isomorphic to the cohomology $H^{p}(\nu_{C,\eps},\nu_{C,\eps}\backslash{}V_{C,\eps};\F)$, which,
via the Thom isomorphism, is isomorphic to the right hand side of \refe{kerlapnce}.
\end{proof}


\section{Proof of the Kirwan-Novikov inequalities} \Label{S:proof}

In this section we prove \reft{cylinder} (and, hence, in view of \refp{cyl=bound}, also
\reft{main}). The main ingredient of the proof is \reft{model} which estimates the spectrum of
the Laplacian $\lap_T$ (cf. \refss{witM}) in terms of the spectrum of the Laplacian
$\bigoplus_{C\in\bfC}\lap_T^\nce$ (cf. \refss{witN}). We postpone the proof of \reft{model} to
the next section.

\subsection{The deformed Laplacian}\Label{SS:addnot}
Set
\begin{gather}\notag
    \tilN= \bigsqcup_{C\in \bfC}\nce, \quad B_r=\bigsqcup_{C\in \bfC} B_{C,r},
    \quad S_r=\bigsqcup_{C\in \bfC} S_{C,r},\\
    g_T^\tilN=\bigoplus_{C\in \bfC}g_T^\nce,\quad h_T^\tilN=\bigoplus_{C\in \bfC}h_T^\nce,
    \quad \lap_T^\tilN= \bigoplus_{C\in\bfC}\lap_T^\nce.\notag
\end{gather}

Recall that the function $\phi(T,y)$ was defined  in \refss{nbE}. Since $\phi(T,y)\equiv T$
outside $B_{C,2}$, we can define  a closed cohomologious to $\ome$ 1-form $\ome_T'\in
\Ome^1(\tilM)$ by the formula
\[
    \ome_T' \ = \ \begin{cases}
     \frac1T\,d\,\big(\,\phi({T},y)\,\tilh_C\,\big)
     \qquad&\text{on}\quad B_{C,5}, \quad \text{for all} \quad C\in \bfC,\\
     \ome \qquad &\text{on} \quad \tilM\backslash B_2.
    \end{cases}
\]
Clearly, the critical set of $\ome'_T$ coincides with $\bfC$.

Consider the one parameter family of connections
\eq{nTprime}
    \n_T' \ = \ \n \ + \ T\,e(\ome_T').
\end{equation}

By construction, the metrics $g_T^{\tilN}$ and $h^{\tilN}_T$ coincide on $B_{6}\setminus
B_{3}\subset\nce\cap{}M$ with $g_1^{\tilN}$ and $h^{\tilN}_1$ respectively. Hence, they can be
extended to metrics $g_T^{\tilM}$, $h_T^{\tilM}$ on the whole $\tilM$ in such a way that the
extensions do not depend on $T$ outside of $B_3$ and the restrictions of these metrics to the
cylindrical part $\calT$ have the form described in \refss{metT}. For each $T\in \RR$ we denote
by $\n_T^{\prime *}$ the formal adjoint of $\n_T'$ with respect to the metrics $g_T^{\tilM}$ and
$h_T^{\tilM}$. Let
\eq{DelTprime}
    \Del_T' \ = \ \frac1T\, \big(\, \n_T'\n_T^{\prime *}+\n_T^{\prime *}\n_T'\, \big).
\end{equation}
Let $\lap^{\prime p}_T, \, \lap^{\tilN,p}_T$ denote the restrictions of operators
$\lap_T',\lap_T^\tilN$ to the space of $p$-forms.

\subsection{Asymptotic of the spectrum of $\lap_T'$} \Label{SS:model}
Let $A$ be a self-adjoint operator with discrete spectrum. For any $\lam>0$, we denote by
$N(\lam,A)$ the number of the eigenvalues of $A$  not exceeding $\lam$ (counting multiplicity).

The following theorem plays a central role in our proof of \reft{cylinder}.
\th{model}Let $\lam_p$ $(p=0,1\nek \dim M)$ be the smallest
   non-zero eigenvalue of $\lap_1^{\tilN,p}$. Then for any $\lam_p>\eps>0$ there exists
   $T_\eps>0$ such that for all $T>T_\eps$
   \begin{equation}\Label{E:DelE=M}\notag
    N(\lam_p-\eps,\lap^{\prime p}_{T})=\dim\Ker\lap_1^{\tilN,p}.
   \end{equation}
\eth
We will prove this theorem in the next section.

\subsection{Proof of \reft{cylinder}} \Label{SS:proof}
By \reft{model},
\begin{equation}\Label{E:N=ker}
    N(\lam_p/2, \lap^{\prime p}_{T})=\dim\Ker \lap^{\tilN,p}_{1},
        \qquad\qquad p=0,1\nek \dim M.
\end{equation}

Let $\tilF_T'$ denote the flat vector bundle $(\tilF,\n_T')$. Since $\ome$ is cohomologious to
$\ome_T'$, the cohomology groups $H^p(\tilM,\tilU_c^-;\tilF_T)$ and $H^p(\tilM,\tilU_c^-;\tilF_T')$
are isomorphic for all $T\in \RR$.

Let us fix a generic $T>0$, cf. \refss{calN}. Then the dimension of the cohomology \/
$H^p(\tilM,\tilU_c^-;\tilF_T')\simeq H^p(\tilM,\tilU_c^-;\tilF_T)$ is equal to the generalized
Novikov number $\bet_i(\xi,\F)$, cf. \refss{calN}. By \refp{tilcomp}, these numbers are equal to
the dimension of the cohomology of the complex
\begin{equation}\Label{E:com-ta}\notag
 \begin{CD}
    0\ \to\ \tilOme^0_T(\tilM,\tilF)@>{\nb'}>>\tilOme^1_T(\tilM,\tilF)@>{\nb'}>>\cdots
            @>{\nb'}>>\tilOme^n_T(\tilM,\tilF)\ \to\ 0.
 \end{CD}
\end{equation}

Let $E^p_T \ (p=0,1\nek \dim M)$ be the subspace of $\tilOme^p_T(\tilM,\tilF)$ spanned by the
eigenvectors of $\lap^{\prime p}_T$ corresponding to the eigenvalues $\lam\le \lam_p/2$. From
\refe{N=ker} and \reft{lapE}(2), we obtain
\begin{equation}\Label{E:k=h}
    \dim E^p_T \ =\
    \dim {H}^{p-\ind(C)}
           (V_{C,\eps},\F_{|_{V_{C,\eps}}}\otimes o(\Sig_C)_{|_{V_{C,\eps}}}).
\end{equation}
Set
\[
    \calM_{\eps}(\lam) \ = \ \sum\, \lam^p\, \dim E^p_T \ = \
    \sum_{C\in \bfC}\, \lam^{\ind{C}}\, \sum\, \lam^p\,
    \dim {H}^{p-\ind(C)}
           (V_{C,\eps},\F_{|_{V_{C,\eps}}}\otimes o(\Sig_C)_{|_{V_{C,\eps}}}).
\]

Since the operator $\lap_T'$ commutes with $\nb'$, the pair $(E^\b_T,\nb)$ is a subcomplex of
\refe{com-ta} and the inclusion induces an isomorphism of cohomology
\[
    H^\b(E^\b_T,\nb') \ \cong \ H^\b(\tilOme^\b(\tilM,\tilF),\nb').
\]
Hence,
\begin{equation}\Label{E:h=b}
    \dim H^p(E^\b_T,\nb') \ =\ \bet_p(\xi,\F), \qquad
            p=0,1\nek \dim M.
\end{equation}
The standard Morse theory arguments  (cf., e.g.,  \cite{Bott2}) imply now that there exists a
polynomial $\calQ_\eps(\lam)$ with non-negative integer coefficients such that
\eq{Morseeps}
    \calM_\eps(\lam) \ - \  \calN_{\ome,\calF}(\lam) \ = \ (1+\lam)\,\calQ_\eps(\lam).
\end{equation}
Using the equality $\bigcap_{r>0}V_{C,\eps}= C$ and the continuity of the \Cech cohomology
(\cite[Ch.~VIII \S6.18]{Dold}), we get:
\[
    \lim_{\eps\to 0}\,
    \dim {H}^{p-\ind(C)}
           (V_{C,\eps},\F_{|_{V_{C,\eps}}}\otimes o(\Sig_C)_{|_{V_{C,\eps}}})
    \ = \ \dim \check{H}^{p-\ind(C)}(C,\F_{|_{C}}\otimes o(C)).
\]
Therefore, $\lim_{\eps\to 0} \calM_\eps(\lam)= \calM_{\ome,\calF}(\lam)$. Letting $\eps\to 0$ in
\refe{Morseeps} we obtain \reft{cylinder}.  \hfill $\square$

\renewcommand{\H}{\text{\( \lap^p_T\)}}
\section{Comparison between the Laplacians on the manifold and on the normal bundle} \Label{S:model}

In this section we prove \reft{model}. Our strategy will be to apply the IMS localization
formula (cf. \cite{CFKS,Shubin96Morse}).

To simplify the notation we will omit the prime and will denote by $\Del_T$ the operator defined in
\refe{DelTprime}.

\subsection{Estimate from above on $N(\lam_p-\eps,\H)$} \Label{SS:above}
We will first show that
\begin{equation}\Label{E:above}
    N(\lam_p-\eps,\Del^p_T)\le\dim\Ker\lap_1^{\tilN,p}.
\end{equation}
by estimating the operator \H\ from below. We will use the technique of \cite{Shubin96Morse},
adding some necessary modifications (see also \cite{BrFar1}).

Recall that the notations $\tilN,\,B_r,\, S_r$, etc. were defined in \refss{addnot}. Using
\refr{C3} we obtain an embedding $i:B_r\to \tilM, \ 0<r<6$.

Let $\tau:\tilN\to \RR$ be a smooth function, such that $\tau(y)\le 1$ for all $y\in B_1$ and
\[
    \tau(x,t)  \ = t \qquad \text{for} \quad (x,t)\in S_{C,1}\times[1,\infty)
    \ \subset \ \tilN.
\]
Then $B_r= \tau^{-1}([0,r))$ for all $r> 1$.

Let us fix a  $C^\infty$ function $j:[0,+\infty)\to [0,1]$ such that $j(s)=1$ \/ for \/  $s\le
4$, $j(s)=0$ \/ for \/  $s\ge 5$ and the function $(1-j^2)^{1/2}$ is $C^\infty$. We define
functions $J,\oJ\in C^\infty(\tilN)$ by
\[
    J(y) \ = \ j(\tau(y)); \qquad \oJ(y) \ = \ \big(1-j(\tau(y))^2\big)^{\frac12}.
\]
Using the diffeomorphism $i:B_6\to i(B_6)$ we can and we will consider $J,\oJ$ as functions on
$\tilM$.

We identify the functions $J,\oJ$ with the corresponding multiplication operators. For operators
$A,B$, we denote by $[A,B]=AB-BA$ their commutator.

The following version of IMS localization formula (cf. \cite{CFKS}) is due to Shubin \cite[Lemma
3.1]{Shubin96Morse}.
\lem{sh1}The following operator identity holds
  \begin{equation}\Label{E:sh1}
    \H=\oJ \H\oJ+J\H J+\frac12[\oJ,[\oJ,\H]]+\frac12[J,[J,\H]].
  \end{equation}
\elem
\begin{proof}
  Using the equality $J^2+\oJ^2=1$ we can write
  $$
    \H=J^2\H+\oJ^2\H=J\H J +\oJ\H\oJ+J[J,\H]+\oJ[\oJ,\H].
  $$
   Similarly,
  $$
    \H=\H J^2+\H\oJ^2=J\H J +\oJ\H\oJ-[J,\H]J-[\oJ,\H]\oJ.
  $$
  Summing these identities and dividing by 2 we come to \refe{sh1}.
\end{proof}
We will now estimate each one of the summands in the right hand side of \refe{sh1}.
\lem{sh3}
  There exist $c>0,\ T_0>0$ such that, for any $T>T_0$,
  \begin{equation}\Label{E:sh3}
    \oJ\H\oJ\ge cT\oJ^2 I.
  \end{equation}
\elem
\begin{proof}
Let $\eta$ be in $L^2\Ome^p((\tilM,\tilF),g_T)$. Using \refe{DelT=} and \refl{ab-homg} we can
write
\eq{oJeta}
    \< \oJ\lap_T^p \oJ\eta,\eta \> = \frac1T\< \lap_T^p \oJ\eta,\eta \> + \< a(x,t) \oJ \eta, \oJ \eta \>+
    T|\ome|^2\oJ^2\<\eta,\eta\>
\end{equation}
(Since the support of $\oJ\eta$ belongs to the set $M\backslash{B_{4}}$, where the scalar
product does not depend on $T$, we can omit this index). The first term of \refe{oJeta} is
positive, the second is bounded, so the last term dominates and the estimate \refe{sh3} holds.
\end{proof}

Let $P^p_T:L^2\Ome^p\big(\tilN,\Ff\big)\to \Ker\Del^{\tilN,p}_T$ be the orthogonal projection.
This is a finite rank operator on $L^2\Ome^p\big(\tilN,\Ff\big)$ and its rank equals
$\dim\Ker\Del^{\tilN,p}_T= \dim\Ker\Del^{\tilN,p}_1$. Clearly,
\begin{equation}\Label{E:ge}
    \Del^{\tilN,p}_T \ + \ \lam_p\, P^p_T \ \ge \ \lam_p I.
\end{equation}
Using the identification $i:B_6\to i(B_6)$ we can consider $JP^p_TJ$  and $J\Del^{\tilN,p}_T J$
as operators on $\Ome^p(\tilM,\tilF)$.

By construction of the operators $\Del_T$ and $\Del_T^{\tilN}$ we have
\[
    J\H J\ =\ J\Del^{\tilN,p}_TJ.
\]
Hence, \refe{ge} implies the following
\lem{local}
  For any $T>0$
  \begin{equation}\Label{E:local}
    J\H J \ + \ \lam_p\,JP^p_TJ \ \ge\  \lam_p\,J^2I, \qquad
        \rk JP^p_TJ\le \dim\Ker\Del^{E,p}_1.
  \end{equation}
\elem
For an operator $A:\Ome^p(\tilM,\tilF)\to \Ome^p(\tilM,\tilF)$, we denote by $\|A\|_T$ its norm
with respect to $L^2$ scalar product on $\Ome^p(\tilM,\tilF)$ defined by metrics $g_T^\tilM$ and
$h_T^\tilM$
\lem{sh2}There exists $C>0$ such that
  \begin{equation}\Label{E:sh2}
    \big\|\,[J,[J,\H]\,\big\|_T \ \le \  CT^{-1}, \qquad\qquad T>0.
  \end{equation}
\elem
\begin{proof}
  The leading symbol of $\H$ is equal to the leading symbol of
  $\frac1{T}\Del$ (here $\Del=\n\n^*+\n^*\n$ is the standard
  Laplacian). Hence (cf., for example, \cite[Proposition~2.3]{BeGeVe})

  $$
    \|[J,[J,\H]\|_T \ = \ -\frac2{T}|dJ|^2 \ \le \ -\, \frac{2}{T}\, \max_{y\in \tilM}\, |dJ(y)|.
  $$
Here we can also omit the index $T$ of the norm of $dJ$, since $\supp dJ\in B_5\setminus B_4$.
\end{proof}

Similarly, one shows that
\begin{equation}\Label{E:sh2oJ}
     \big\|\,[\oJ,[\oJ,\H]\,\big\|_T\ \le\ CT^{-1}.
\end{equation}
{}From \refl{sh1}, \refl{sh3}, \refl{local}, \refl{sh2} and \refe{sh2oJ} we get the following
\cor{above}
   For any $\eps>0$, there exists $T_0>0$ such that for any $T>T_0$
   \begin{equation}\Label{E:a'}
    \H\ +\ \lam_p\,JP^p_TJ \ \ge \  (\lam_p-\eps)\,I, \qquad
        \rk JP^p_TJ \ \le\ \dim\Ker\Del^{E,p}_T.
   \end{equation}
\ecor
The estimate \refe{above} follows now from \refc{above} and the following general lemma \cite[p.
270]{ReSi4}.
\lem{general}
   Assume that $A, B$ are self-adjoint operators in a Hilbert space
   $\calH$ such that  $\rk B\le k$ and there exists $\mu>0$
   such that
   $$
    \langle (A+B)u,u\rangle \ge \mu\langle u,u\rangle
        \quad \text{for any} \quad u\in\Dom(A).
   $$
   Then $N(\mu -\eps, A)\le k$ for any $\eps>0$.
\elem

\subsection{Estimate from below on $N(\lam_p-\eps, \Del^p_T)$} \Label{SS:below}
To prove \reft{model} it remains to  show  that
\begin{equation}\Label{E:below}
    N(\lam_p-\eps,\H) \ \ge\ \dim\Ker\Del^{\tilN,p}_T.
\end{equation}
Let $E^p_T$ be the subspace of $\Ome^p(\tilM,\tilF)$ spanned by the eigenvectors of $\Del^p_T$
corresponding to the eigenvalues $\lam\le \lam_p-\eps$ and let $\Pi^p_T:\Ome^p(\tilM,\tilF)\to
E_T^p$ be the orthogonal projection. Then
\begin{equation}\Label{E:rkPi}
    \rk \Pi^p_T \ =\ N(\lam_p-\eps,\Del^p_T).
\end{equation}
Using the diffeomorphism $i:B_6\to i(B_6)$ we can consider $J\Pi^p_TJ$ as an operator on
$L^2\Ome^p\big(E,\Ff\big)$. The proof of the following lemma does not differ from the proof of
\refc{above}.
\lem{below} For any $\del>\eps$, there exists $T>0$ such that for any
   $t>T$
   \begin{equation}\Label{E:a''}
    \Del^{E,p}_T\ +\ J\Pi^p_TJ\ \ge\ (\lam_p-\del)\,I.
   \end{equation}
\elem
The estimate \refe{below} follows now from \refe{rkPi}, \refl{below} and \refl{general}. \hfill
$\square$


\providecommand{\bysame}{\leavevmode\hbox to3em{\hrulefill}\thinspace}
\providecommand{\MR}{\relax\ifhmode\unskip\space\fi MR }
\providecommand{\MRhref}[2]{%
  \href{http://www.ams.org/mathscinet-getitem?mr=#1}{#2}
} \providecommand{\href}[2]{#2}


\begin{thebibliography}{10}

\bibitem{BrAlesker}
S.~Alesker and M.~Braverman, \emph{{Cohomology of a Hamiltonian T-space with
  involution}}, Preprint.

\bibitem{AtBott82}
M.~F. Atiyah and R.~Bott, \emph{The {Yang}-{Mills} equations over {Riemann}
  surface}, Phil. Trans. R. Soc. London, ser A \textbf{308} (1982), 523--615.

\bibitem{AtBott84}
\bysame, \emph{The moment map and equivariant cohomology}, Topology \textbf{23}
  (1984), 1--28.

\bibitem{BeGeVe}
N.~Berline, E.~Getzler, and M.~Vergne, \emph{Heat kernels and {Dirac}
  operators}, Springer-Verlag, 1992.

\bibitem{Bis1}
J.-M. Bismut, \emph{The {Witten} complex and the degenerate {Morse}
  inequalities}, J. Dif. Geom. \textbf{23} (1986), 207--240.

\bibitem{Bott2}
R.~Bott, \emph{Morse theory indomitable}, Publ. Math. IHES \textbf{68} (1988),
  99--114.

\bibitem{BrFar1}
M.~Braverman and M.~Farber, \emph{The {Novikov}-{Bott} inequalities}, C.R.
  Acad. Sci. Paris \textbf{t.~321, S\'erie I} (1995), 897--902.

\bibitem{BrFar3}
\bysame, \emph{Novikov inequalities with symmetry}, C. R. Acad. Sci. Paris
  S\'er. I Math. \textbf{323} (1996), 793--798.

\bibitem{BrFar4}
M.~Braverman and M.~Farber, \emph{Equivariant {N}ovikov inequalities},
  $K$-Theory \textbf{12} (1997), 293--318.

\bibitem{BrFar2}
\bysame, \emph{Novikov type inequalities for differential forms with
  non-isolated zeros}, Math. Proc. Cambridge Philos. Soc. \textbf{122} (1997),
  357--375.

\bibitem{BrMiSh02}
M.~Braverman, O.~Milatovich, and M.~Shubin, \emph{Essential selfadjointness of
  {S}chr\"odinger-type operators on manifolds}, Russian Math. Surveys
  \textbf{57} (2002), 41--692.

\bibitem{Chernoff73}
P.~Chernoff, \emph{Essential self-adjointness of powers of generators of
  hyperbolic equations}, J. Functional Analysis \textbf{12} (1973), 401--414.

\bibitem{Conley78}
C.~Conley, \emph{Isolated invariant sets and the {M}orse index}, CBMS Regional
  Conference Series in Mathematics, vol.~38, American Mathematical Society,
  Providence, R.I., 1978.

\bibitem{CFKS}
H.L. Cycon, R.G. Froese, W.~Kirsch, and B.~Simon, \emph{Schr\"odinger operators
  with applications to quantum mechanics and global geometry}, Texts and
  Monographs in Physics, Springer-Verlag, 1987.

\bibitem{Dold}
A.~Dold, \emph{Lectures on algebraic topology}, Classics in Mathematics,
  Springer-Verlag, Berlin, 1995, Reprint of the 1972 edition.

\bibitem{EliashGromov98}
Y.~Eliashberg and M.~Gromov, \emph{Lagrangian intersection theory:
  finite-dimensional approach}, Geometry of differential equations, Amer. Math.
  Soc. Transl. Ser. 2, vol. 186, Amer. Math. Soc., Providence, RI, 1998,
  pp.~27--118.

\bibitem{FarberShustin00}
M.~Farber and E.~Shustin, \emph{Witten deformation and polynomial differential
  forms}, Geom. Dedicata \textbf{80} (2000), no.~1-3, 125--155.

\bibitem{Farber85}
M.S. Farber, \emph{Exactness of the {Novikov} inequalities}, Functional Anal.
  Appl. \textbf{19} (1985), 40--48.

\bibitem{Floer89}
A.~Floer, \emph{Witten's complex and infinite dimensional {Morse} theory}, J.
  Diff. Geom \textbf{30} (1989), 207--221.

\bibitem{GromLaw83}
M.~Gromov and B.~Lawson, \emph{Positive scalar curvature and the {D}irac
  operator on complete {R}iemannian manifolds}, Inst. Hautes \'Etudes Sci.
  Publ. Math. (1983), no.~58, 295--408.

\bibitem{Hirsch}
M.~Hirsch, \emph{Differential topology}, Graduate Texts in Math., vol.~33,
  Springer, Berlin, 1976.

\bibitem{Kirwan84}
F.~C. Kirwan, \emph{Cohomology of quotients in symplectic and algebraic
  geometry}, Mathematical Notes, vol.~31, Princeton University Press,
  Princeton, NJ, 1984.

\bibitem{Novikov81}
S.P. Novikov, \emph{Multivalued functions and functionals. {An} analogue of the
  {Morse} theory}, Soviet Math. Dokl. \textbf{24} (1981), 222--226.

\bibitem{Novikov82}
\bysame, \emph{The {Hamiltonian} formalism and a multivalued analogue of
  {Morse} theory}, Russian Math. Surveys \textbf{37} (1982), 1--56.

\bibitem{Pazhitnov87}
A.~Pazhitnov, \emph{An analytic proof of the real part of {Novikov's}
  inequalities}, Soviet Math. Dokl. \textbf{35} (1987), 456--457.

\bibitem{ReSi4}
M.~Reed and B.~Simon, \emph{Methods of modern mathematical physics {I}{V}:
  {Analysis} of operators}, Academic Press, London, 1978.

\bibitem{Shubin96Morse}
M.~A. Shubin, \emph{Semiclassical asymptotics on covering manifolds and {M}orse
  inequalities}, Geom. Funct. Anal. \textbf{6} (1996), 370--409.

\bibitem{Shubin99}
\bysame, \emph{Spectral theory of the {S}chr\"odinger operators on non-compact
  manifolds: qualitative results}, Spectral theory and geometry (Edinburgh,
  1998), Cambridge Univ. Press, Cambridge, 1999, pp.~226--283.

\bibitem{Viterbo92}
C.~Viterbo, \emph{Symplectic topology as the geometry of generating functions},
  Math. Ann. \textbf{292} (1992), no.~4, 685--710.

\end{thebibliography}
\end{document}